\documentclass[a4paper]{article}
\usepackage[latin1]{inputenc}  \usepackage[T1]{fontenc}
\usepackage{lmodern}           

\usepackage[italian,spanish,german,frenchb,english]{babel}
\usepackage{theorem}
\usepackage{amssymb}
\usepackage{amsfonts}
\usepackage{amsmath}
\usepackage{amsxtra}

\usepackage{pstricks}
\usepackage{showidx}
\usepackage{amscd}
\usepackage[active]{srcltx}
\usepackage{multicol}
\usepackage{fancyhdr}
\usepackage{changebar}
{\theoremstyle{change} \theoremheaderfont{\normalfont\bfseries}
\theorembodyfont{\slshape}
\newtheorem{Prop}{Proposition:}[section]}

{\theoremstyle{change} \theorembodyfont{\slshape}
\newtheorem{Theo}[Prop]{Theorem:}}

{\theoremstyle{change} \theorembodyfont{\slshape}
\newtheorem{Cor}[Prop]{Corollary:}}

{\theoremstyle{change} \theorembodyfont{\slshape}
\newtheorem{Lem}[Prop]{Lemma:}}

{\theoremstyle{change} \theorembodyfont{\upshape}
\newtheorem{Rem}[Prop]{\normalfont\scshape {Remark:}}}

{\theoremstyle{change} \theorembodyfont{\upshape}
\newtheorem{Defi}[Prop]{\normalfont\scshape{Definition:}}}

{\theoremstyle{change} \theorembodyfont{\upshape}
}

{\theoremstyle{change} \theorembodyfont{\upshape}
\newtheorem{Not}[Prop]{\normalfont\scshape{Notation:}}}

{\theoremstyle{change} \theorembodyfont{\upshape}
\newtheorem{Nr}[Prop]{\kern -1ex}}

{\theoremstyle{change} \theorembodyfont{\upshape}
}

\newcommand{\qed}{\hfill \mbox{\raggedright \rule{.07in}{.1in}}}

\catcode`\á=\active \def á{\'a}
 \catcode`\Á=\active \def Á{\'A}
  \catcode`\ó=\active \def ó{\'o}
  \catcode`\é=\active \def é{\'e}
  \catcode`\ú=\active \def ú{\'u}
  \catcode`\í=\active \def í{\'{\i}}
  \catcode`\ñ=\active \def ñ{\~n}
  \catcode`\Ñ=\active \def Ñ{\~N}
  \catcode`\¿=\active \def ¿{?`}
  \catcode`\º=\active \def º{$^{\underline{o}}$}
  \catcode`\ª=\active \def ª{$^{\underline{a}}$}
  \catcode`\¡=\active \def ¡{!`}
  \catcode`\â=\active \def â{\^{a}}
  \catcode`\ê=\active \def ê{\^{e}}
  \catcode`\î=\active \def î{\^{\i}}
  \catcode`\ô=\active \def ô{\^{o}}
  \catcode`\û=\active \def û{\^{u}}
  \catcode`\ç=\active \def ç{\c{c}}
  \catcode`\ü=\active \def ü{\"{u}}
  \catcode`\ö=\active \def ö{\"{o}}

 \newcommand{\ee}{\mathbb{E}}

 \newenvironment{dem}{\noindent{\it Proof}.}{\qed}

\newcommand{\oei}{\overset{\circ}{E_i}}

\newcommand{\pro}{\mathbb{P}}
\newcommand{\oo}{\mathcal{O}}

\newcommand{\mm}{\mathfrak{m}}

\newcommand{\Nu}{{ \mbox{{\LARGE $\nu$}} } }

\title{Generalised Poincar\'e series and \\ embedded resolution of curves}
\author{Julio Jos\'e Moyano-Fern\'andez    \thanks{Supported partially by the grant of the Spanish Government
Ministerio de Educaci\'on MTM2007--64704, in cooperation with The European Union in the framework of the founds  ``FEDER'', by the grant of the
regional Government Junta de Castilla y Le\'on VA065A07, by
the grant of the Deutsches Akademischen Austauschdienst (DAAD)--La
Caixa, and by the Deutsche Forschungsgemeinschaft (DFG). The author is thankful to Prof. Dr. Karlheinz Kiyek  and
Prof. Dr. F\'elix Delgado for nice conversations and useful remarks.
He is also thankful to the Universities of Paderborn and
Valladolid for kind hospitality.}}

\date{Institut f\"ur Mathematik, Universit\"at
Osnabr\"uck\\
Email:  jmoyanof@uni-osnabrueck.de}
\begin{document}
\maketitle

\renewcommand{\abstractname}{Abstract}
\begin{abstract}
 The purpose of this paper is to extend the notions of generalised
Poincar\'e series and divisorial generalised Poincar\'e series (of motivic nature)
introduced by Campillo, Delgado and Gusein--Zade for complex curve singularities to curves defined
over perfect fields, as well as to express them in terms of an
embedded resolution of curves.
\end{abstract}

\noindent
AMS-Classification: 14H20, 32S99

\section{Introduction}

A. Campillo, F. Delgado and K. Kiyek introduced in 1994 a
multivariable Poincar\'e series $P(t_1, \ldots , t_r)$ (from now on
denoted by $P(\underline{t})$) associated with the valuations of
the integral closure of a one--dimensional local Cohen--Macaulay
ring (cf. \cite[(3.8)]{cadeki}). For the case of valuations of a
complex plane curve singularity, Campillo, Delgado and Gusein--Zade
interpreted the Poincar\'e series $P(\underline{t})$ as an integral
over the local ring of germs with respect to the Euler
characteristic, expressing it also in terms of an embedded
resolution of curves (\cite{cadegu4}, \cite{cadegu6},
\cite{duke}). This new approach allowed them also to define the
Poincar\'e series $\widehat{P}(\underline{t})$ associated with the
extended semigroup of the curve singularity just by taking the integral
over the extended semigroup (cf. \cite{duke}). The same philosophy
can be applied to divisorial valuations to get a divisorial
Poincar\'e series $P^D(\underline{t})$ of the divisorial value
semigroup (i.e., the value semigroup arising by considering
divisorial valuations), or even a semigroup divisorial Poincar\'e
series $\widehat{P}^D (\underline{t})$ if we take the divisorial
extended semigroup of the singularity (see \cite{Delgado4}, \cite{Delgado5}).
\medskip

We can also take the generalised Euler characteristic instead of
the classical one; it leads to the study of the motivic versions
of the previous series, namely: the generalised Poincar\'e series
$P_g(\underline{t})$; the generalised semigroup Poincar\'e series
$\widehat{P}_g (\underline{t})$; the generalised divisorial
Poincar\'e series $P^D_g(\underline{t})$; and the generalised
divisorial semigroup Poincar\'e series
$\widehat{P}^D_g(\underline{t})$. They all were studied in
\cite{cadegu11} for complex curve singularities.
\medskip

An interesting question is how to translate these
ideas in more general contexts than $\mathbb{C}$. Some work in
this direction has been already done (see \cite{demo} for Poincar\'e
series over finite fields, and \cite{kimo} for
$P^D(\underline{t})$ and $P^D_g (\underline{t})$ over non--finite
fields). Also, a related motivic zeta function was treated by W. Z\'u\~niga and the author in \cite{mozu}.
Following this direction, we wish to investigate in the present paper the generalised Poincar\'e
series $P_g(\underline{t})$, the generalised divisorial Poincar\'e
series $P^D_g(\underline{t})$ and the generalised divisorial
semigroup Poincar\'e series $\widehat{P}^D_g(\underline{t})$
attached to curves defined over perfect fields; in particular we want to
express them---in this more general context---in terms of an
embedded resolution of the curve. In the process we include some
omissions in \cite{cadegu11}, which we will opportunely indicate.
As for prerequisites, the reader is expected to be familiar with
the theory of two--dimensional regular local rings, especially with
the concept of ideal transform. For a treatment of this topic we
refer the reader to \cite{kiyek}.
\medskip

The paper is organised as follows. Sect.~\ref{sec:2} establishes
the terminology and contains a brief summary about embedded
resolution of curves. Sect. \ref{sec:int} deals with the
Grothendieck ring of the category of quasi--projective schemes of
finite type over a perfect field $k$, the generalised Euler
characteristic and the definition of generalised Poincar\'e series,
according to the exposition presented in \cite{demo}. The aim of
Sect.~\ref{section:gps} is to give a formula expressing the
generalised Poincar\'e series in terms of an embedded resolution of
curves. First, we define a semigroup homomorphism (the map
\textsf{Init} of Subsect.~\ref{subsection:init}) which will be a
locally trivial fibration over each connected component of its
image. We study also the preimage
$\mathrm{\textsf{Init}}^{-1}(\varsigma)$ for every element
$\varsigma$ on the connected component of the image (Corollary
\ref{cor:F(n)}), and compute its (finite) codimension by means of
the Hoskin--Deligne formula (see Subsect.~\ref{subsection:codimension}; 
the explicit calculation is provided
in Proposition \ref{prop:efedeene}). Using all these data and
applying Fubini's formula to the map \textsf{Init} we get a
description of the generalised Poincar\'e series in terms of the
embedded resolution (Theorem \ref{thm:uno}). Finally, in Sect.~\ref{sec:div} 
we introduce both the generalised divisorial
Poincar\'e series and the generalised divisorial semigroup Poincar\'e
series of a curve defined over a perfect field and proceed with a
similar study; the first series is described in Theorem
\ref{thm:gendiv1}, whereas the second one needs some extra work: we show
in \ref{nr:56} the analogous of the map \textsf{Init} in this
context; it is in fact an isomorphism (Lemma \ref{lemma:iso}). The
final formula is given by Theorem \ref{thm:divisorialintermsof}.
In the rest of the section we explain briefly some conventions,
terminology and basic definitions to be used along this paper.

\section{Preliminaries on blowing--ups and exceptional divisors} \label{sec:2}

\Nr Let $R$ be a two--dimensional regular local ring having a
maximal ideal $\mathfrak{m}$ and a perfect residue field $k_R$. Let $f \in R \setminus \{0 \}$
be a non--unit reduced element. We say
that $f$ defines a \textbf{curve} $C_f$ (or simply $C$ when no confusion can arise) on $R$. The irreducible
factors of $f$, let us say $f=f_1 \cdot \ldots \cdot f_r$, are
called the \emph{irreducible components} of $C$. Sometimes we will
refer to this fact expressing $C=C_1 \cup \ldots \cup C_r$. The
curve is \textbf{reduced} (resp. \textbf{analytically reduced}) if
$R/fR$ is reduced (resp., if $\widehat{R}/f\widehat{R}$ is
reduced, where $\widehat{R}$ stands for the $\mm$-adic completion
of $R$).
\medskip

\Nr \label{Manis} Let $C=C_1 \cup \ldots \cup  C_r$ be an
analytically reduced curve defined by the element $f=f_1 \cdot \ldots \cdot
f_r \in R \setminus \{ 0 \}$. Take the ring $R/fR$ and its
integral closure $\overline{R/fR}$. It defines finitely many
discrete Manis valuations $v_i$ having $\overline{R/f_i R}$ as
discrete Manis valuation ring, for all $1 \le i \le r$ (cf.
\cite[Chap. I, (2.2)]{kiyek}). Note that the $v_1, \ldots , v_r$
are not valuations of the ring $R$ (the preimage of $\infty$ via
$v_i$ is $f_i R$, for $1 \le i \le r$), but they allow us to
define a multi--index filtration of the ring $R$ given by the
ideals $J(\underline{n}):=\{ z \in R \mid v_i (z) \ge n_i, ~ ~ 1
\le i \le r \}$, for $\underline{n}:= (n_1 \ldots, n_r) \in
\mathbb{Z}^r$.
\medskip

\Nr The multiplicity intersection of two curves $C_f,C_g$ given by
elements $f,g \in R \setminus \{ 0\}$ is just the multiplicity
intersection of these elements:
\[
(C_f \cdot C_g)=(f \cdot g)_R:= \ell_R (R/(f,g)).
\]
\medskip

\Nr A curve $C$ is said to be normal crossing if it is regular, or
if $C$ has two regular components $C_1, C_2$, defined by two elements $f_1,f_2
\in  R \setminus \{ 0\}$ resp., such that $\ell_R \left ( R / (f_1 R, f_2 R) \right
) = 1$. In this case, $\{f_1,f_2 \}$ is a regular system of
parameters of $R$.
\medskip

\Nr \label{seq} Let $X_0:=\mathrm{Spec}(R)$ be a regular scheme of finite type over $k$ of dimension two; take a closed point $p_0 \in X_0$ and blow up at $p_0$ to get another two--dimensional regular scheme $X_1$. By repeating the process $s$ times we obtain in this way a finite sequence of blowing-ups 
\[
X:=X_{s} \overset{\pi_{s}}{\longrightarrow} X_{s-1}
\overset{\pi_{s-1}}{\longrightarrow} \ldots \longrightarrow X_3
\overset{\pi_3}{\longrightarrow} X_2
\overset{\pi_2}{\longrightarrow} X_1
\overset{\pi_1}{\longrightarrow} X_0,
\]
where $\pi_i$ is the blow up at a closed point $p_{i-1} \in X_{i-1}$, for every $1 \le i \le s$.
Write $\pi := \pi_{s} \circ \pi_{s-1} \circ \ldots \circ \pi_2
\circ \pi_1$. The \textbf{exceptional divisor of $\pi$} is defined to be the
reduced inverse image $\pi^{-1}(\{\mm \})$ of the maximal ideal
$\mm$ of the ring $R$. It coincides with the union of the strict
transforms of the exceptional divisors of each $\pi_j$, $1 \le j
\le s-1$, together with $\pi_{s}^{-1}(\{p_{s-1} \})$, i.e., it
has $s$ different irreducible components, every two components
meet transversally at one point, and no three components meet at a
point. Each irreducible component of $E$ is isomorphic to a scheme
$\mathrm{Proj}(k_j [\overline{x},\overline{y}])$, $1 \le j \le s$,
where $k_j$ is a finite extension of $k_R$ of degree $h_j$ (cf. \cite[Chap.~VII]{kiyek}, \cite[Sect.~2]{moyano2})).
\medskip

\Nr \label{nr:215} We will write $E_{i,i}$ for the exceptional
divisor of $\pi_i$ as divisor of $X_i$, and we denote by $E_{i,j}$
(resp. ${E_{i,j}^{\ast}}$) the strict transform (resp. the total
transform) of $E_{i,i}$ in $X_j$ by the morphism $X_j
\longrightarrow X_i$, for $j>i$. We denote by $E_i$ (resp.
$E_i^{\ast}$) the strict (resp. total) transform $E_{i,s}$ (resp.
${E_{i,s}^{\ast}}$) by the morphism $X \longrightarrow X_i$. 
Let $\ee$ be the subgroup of $1$--cycles of $X$ of the form
$\sum_{i=1}^{s} n_i E_i$, with $n_i \in \mathbb{Z}$. Both the set
$\{E_i \}$ and $\{E_i^{\ast}\}$ are basis of $\ee$.  One has also a symmetric
bilinear intersection form   $\mathbb{E} \times \mathbb{E}  \rightarrow  \mathbb{Z} $ given by intersecting cycles 
$ (A,B)  \mapsto  \mathrm{deg}_X (A \cdot B)$. From the
projection formula follows
\[
\deg_X (E_i^{\ast} \cdot E_j^{\ast}) = -\delta_{ij} h_i,
\eqno(\ddag)
\]
where $\delta_{ij}$ is the Kronecker's delta  (see \cite[Theorem 9.2.12, p.~398]{liu}).
\medskip

\Nr \label{nr:matrices} Let $p_i \in X_i, p_j \in X_j$ with
$(\pi_j \circ \pi_{j-1} \circ \ldots \circ \pi_i)(p_j)=p_i$, for
$j>i$. We say that $p_j$ is proximate to $p_i$, and denote it by
$p_j > p_i$, if and only if we have $p_j \in E_{i,j-1}$
(set-theoretically). The basis change matrix from $\{E_i^{\ast}\}$
to $\{E_i\}$ is called the proximity matrix associated to $\pi$.
This matrix will be denoted by $P_{\pi}$. It is easy to check that
any entry $p_{ij}$ of $P_{\pi}$ is equal to $1$ if $i=j$, to $-1$
if $j>i$ and $0$ otherwise. The matrix of the intersection form
given in \ref{nr:215} in the basis $\{E_i^{\ast} \}$ is
$-\Delta_{\pi}$, being $\Delta_{\pi}$ the $s \times s$--diagonal
matrix with entries the extension degrees $h_1, \ldots, h_s$. Moreover, the matrix of
the intersection form in the basis $\{E_i \}$ is
$N_{\pi}:= -P_{\pi} \cdot \Delta_{\pi} \cdot P_{\pi}^t$,
where $P_{\pi}^t$ denotes the transpose of $P_{\pi}$. The matrix
$N_{\pi}$ is called the intersection matrix (with respect to the
basis $\{E_i \}$) associated with $\pi$.  We will write $P$, $N$
and $\Delta$ instead of $P_{\pi}$, $N_{\pi}$ and $\Delta_{\pi}$
whenever the blow--up given by $\pi$ is clear from the context (cf. \cite[Sect.~4]{moyano2}).
\medskip

\section{Integrals with respect to the generalised Euler
characteristic}\label{sec:int}

\Nr Let $\Nu_{k}$ be the
category of quasi--projective schemes of finite type over a perfect field $k$. The
\textbf{Grothendieck ring} of $\Nu_{k}$, denoted by $K_0
(\Nu_{k})$, is defined to be the free Abelian group on isomorphism
classes $[X]$ of quasi--projective schemes $X$ of finite type over
$k$ subject to the relations (i) $[X_1]=[X_2]$ if $X_1 \cong X_2$ for $X_1,X_2 \in
\Nu_{k}$; (ii) $[X]=[X \setminus Z]+[Z]$ for a closed
subscheme $Z$ of $X \in \Nu_k$; and taking the fibred product as multiplication:
(iii) $[X_1] \cdot [X_2]=[X_1 \times_k X_2]$ for $X_1,X_2 \in
\Nu_k$. Notice that, if $X_1$ and $X_2$ are reduced, then $X_1 \times_k
X_2$ is also reduced, because $k$ is perfect. The neutral element of $K_0 (\Nu_k)$ with respect to the
addition will be denoted by $0$ and corresponds to the class
$[\varnothing]$ of the empty set. The Grothendieck ring $K_0
(\Nu_{k})$ is commutative and with unit, the unit being the class
of $\mathrm{Spec}(k)$. Let $k [T]$ be the polynomial ring in one indeterminate $T$
over the field $k$. We denote $\mathbb{A}^1_k:= \mathrm{Spec}(k[T])$ and
$\mathbb{L}:=[\mathbb{A}^1_k]$, which is a non--zero divisor of $K_0
(\Nu_k)$ (see \cite{poonen}); thus we can consider $\{
\mathbb{L}^n \}_{n \in \mathbb{N}}$ as a multiplicatively closed
subset of $K_0 (\Nu_k)$, and localise with respect to
$\mathbb{L}$ to define $\mathcal{M}_{k}:= K_0
(\Nu_{k})_{\mathbb{L}}$.
\medskip

\Nr We can consider locally closed subspaces $Y$ of $X \in
\Nu_{k}$ as schemes and therefore as elements of $\Nu_{k}$ (cf.
\cite[Proposition 4.6.1., p.~273]{ega1}). Hence the locally closed
subsets have a natural image on the Grothendieck ring $K_0
(\Nu_{k})$. We would like also to be able to speak about classes
of constructible subsets (i.e., finite disjoint unions of locally
closed subspaces) of elements of $\Nu_k$. Then we recall the
following result (see \cite[Introduction]{denef}):
\medskip

\begin{Prop}
If $Y$ is a scheme of finite type over $k$, then the map
$Y^{\prime} \to [Y^{\prime}]$ from the set of closed subschemes of
$Y$ to $K_0(\Nu_k)$ extends uniquely to a map $Z \to [Z]$ from the
set of constructible subsets of $Y$ to $K_0(\Nu_k)$ satisfying
$[Z \cup Z^{\prime}] = [Z]+[Z^{\prime}]-[Z \cap Z^{\prime}]$ 
for $Z, Z^{\prime} \in \Nu_k$.
\end{Prop}
\medskip

Therefore, every constructible subset $Z$ of an element of $\Nu_k$
has a well--defined image $[Z] \in K_0(\Nu_k)$. Note the importance
of constructible subsets of a scheme $X$: They are the smallest
algebra of sets containing the closed sets for the Zariski
topology. Moreover, the Grothendieck ring behaves also well with
respect to trivial fibrations (see \cite[p.~6]{blickle}):
\medskip

\begin{Lem} \label{lem:36}
Let $f:Y \to X$ be a piecewise trivial fibration with constant
fibre $Z$. This means that one can write $X = \bigsqcup X_i$ as a
finite disjoint union of locally closed subsets $X_i$ such that
over each $X_i$ one has $f^{-1} X_i \cong X_i \times Z$ and $f$ is
given by the projection onto $X_i$. Then in $K_0 (\Nu_k)$ we have
$[Y]=[X] \cdot [Z]$.
\end{Lem}
\medskip

\Nr \label{nr:38}The morphism $\phi: \mathbb{A}^1_k \setminus \{0 \}
\times \mathbb{A}^1_k \setminus \{0 \} \to \mathbb{A}^1_k
\setminus \{0 \}$ so that $\phi(a,b)=a \cdot b$ for all $a, b \in
\mathbb{A}^1_k \setminus \{0 \}$ defines a group scheme which
is called the \emph{multiplicative group} and is denoted by
$\mathbb{G}_m$ (see \cite[p.~324]{hartshorne}). We define the
class $[k^{\ast}]$ of the group of units $k^{\ast}$ of $k$ as the
class $[\mathbb{G}_m] \in K_{0}(\Nu_k)$. Since $\mathbb{G}_m \cong
\mathrm{Spec} \left ( k \left [ x, \frac{1}{x} \right ] \right )$,
we get $[\mathbb{G}_m] = \left [ \mathrm{Spec} \left ( k \left [ x
, \frac{1}{x} \right ] \right ) \right ]$, and therefore $[k^{\ast}] = \mathbb{L}-1$.
\medskip

\Nr Let $C$ be a curve with equation $f$ and $\mathcal{O}:=R/fR$.
Let $p$ be a non-negative integer and let $J_{\oo}^p$ be the space
of $p$--jets over $\oo$, which is a $k_R$--vector
space of finite dimension $d(p)$. Let us consider its projectivisation
$\pro J_{\oo}^p$ and let us adjoin one point to this (that is,
$\pro^{\ast}J_{\oo}^p = \pro J_{\oo}^p \cup \{\ast \}$ with $\ast$
representing the added point) in order to have a well--defined map
$\pi_p: \pro \oo \to \pro^{\ast} J_{\oo}^p$. A subset $X \subset
\pro \oo$ is said to be cylindric if there exists a constructible
subset $Y \subset \pro J_{\oo}^p \subset \pro^{\ast} J_{\oo}^p$
such that $X = \pi_p^{-1} (Y)$.
\medskip

The generalised Euler characteristic $\chi_g (X)$ of a cylindric
subset $X$ is the element $[Y] \cdot \mathbb{L}^{-d(p)}$ in the
ring $\mathcal{M}_{k_R}$, where $Y=\pi^{-1}_p(X)$ is a constructible
subset of $\mathbb{P}\oo$. Note that $\chi_g (X)$ is well--defined,
because if $X=\pi^{-1}_q(Y^{\prime})$, $Y^{\prime} \subset \pro
J_{\oo}^q$ and $p \ge q$, then $Y$ is a locally trivial fibration
over $Y^{\prime}$ and therefore $[Y] = [Y^{\prime}] \cdot
\mathbb{L}^{d(p)-d(q)}$.
\medskip

Let $\psi: \mathbb{P} \oo \to G$ be a function with values in an
Abelian group $G$ with countably many values. It is said
\emph{cylindric} if, for each $a \in G \setminus \{ 0 \}$, the set
$\psi^{-1}(a) \subseteq \mathbb{P}\oo$ is cylindric. As it is
defined in \cite{cadegu3}, the \emph{integral} of a cylindric
function $\psi$ over $\mathbb{P}\oo$ with respect to the
generalised Euler characteristic is
\[
\int_{\mathbb{P}\oo} \psi d \chi_g := \sum_{a \in G \setminus \{0
\}} \chi_g \left ( \psi^{-1} (a) \right ) \cdot a,
\]
if this sum makes sense in $\mathcal{M}_{k_R} \otimes_{\mathbb{Z}} G$;
in such a case, the function $\psi$ is said to be
\emph{integrable}.
\medskip

\begin{Defi} \label{poincareseries}
The generalised Poincar\'e series of the multi--index filtration
given by the ideals $J(\underline{n})$ (see \ref{Manis}) is the
integral
\begin{equation} \label{eqn:poincareseriesdef}
P_g (t_1, \ldots ,t_r;\mathbb{L}):=\int_{\mathbb{P} R}
\underline{t}^{\underline{v}(h)} d \chi_g \in \mathcal{M}_{k_R}
\nonumber
\end{equation}
where $\underline{t}^{\underline{v}(h)}:=t_1^{v_1 (h)} \cdot
\ldots \cdot t_r^{v_r(h)}$ is considered as a function on
$\mathbb{P}R$ with values in $\mathbb{Z}[\![t_1, \ldots , t_r]\!]$
(the vector $\underline{t}^{\underline{v}(h)}$ is supposed to be
$\underline{0}$ as soon as at least one of the $v_i(h)$ is
$\infty$). Notice that the sub--index $g$ in $P_g(\underline{t})$
is just notation: It refers to integration with respect to the
\emph{generalised} Euler characteristic.
\end{Defi}

\section{Generalised Poincar\'e series in terms of an embedded
resolution}\label{section:gps}

From now on, we assume that the residue field $k_R$ is isomorphic
to a perfect field $K$ contained in $R$, that is, we assume the
existence of a perfect coefficient field $K$ (it will be only
needed to define the integral over $\pro R$).
\medskip

Let $C$ be a curve with equation $f$ satisfying the same
hypothesis as in \ref{Manis}. Since the curve is analytically
reduced, the total transform of $C$ is a normal crossing curve at
some point of a sequence of quadratic transforms (cf.
\cite[Chap. VII, (8.13), p.~301]{kiyek}). Consequently, we get a finite sequence of blowing-ups
\[
X=X_{s} \overset{\pi_{s}}{\longrightarrow} X_{s-1}
\overset{\pi_{s-1}}{\longrightarrow} \ldots \longrightarrow X_3
\overset{\pi_3}{\longrightarrow} X_2
\overset{\pi_2}{\longrightarrow} X_1
\overset{\pi_1}{\longrightarrow} X_0=\mathrm{Spec}(R)
\]
with $\pi = \pi_{s} \circ \pi_{s-1} \circ \ldots \circ \pi_2
\circ \pi_1$, as in \ref{seq}. Let $E$ be the exceptional divisor of $\pi$. For each $1 \le i \le
s$, let $\overset{\circ}{E_i}$ be the component $E_i$ of the
exceptional divisor of $\pi$ minus the intersection points with
all other components of the total transform of the curve. In the
same way, we consider the set $\overset{\bullet}{E_i}$, which is
just the component $E_i$ of the exceptional divisor minus the
intersection points with other components $E_j$ with $j \ne i$.
Note that both $\overset{\circ}{E_i}$ and $\overset{\bullet}{E_i}$
are quasi--projective schemes.
\medskip

The goal in writing this section is to describe the integral
defining the ge\-neralised Poincar\'e series (cf. Definition
\ref{poincareseries}) in terms of an embedded resolution of
curves. First of all, we will distinguish three different types of
points on the exceptional divisor $E$: The intersection points
among components of $E$, the intersection points between
components of $E$ and the components of the strict transform of
the curve, and the points in $\overset{\circ}{E_i}$. By
considering separately the set of the previous first two kind of
points, as well as the symmetric product of the smooth components
$\overset{\circ}{E_i}$, we will define in Subsect.
\ref{subsection:init} a semigroup $Y$ and a surjection
\[
\mathrm{\textsf{Init}}: \mathbb{P}R^{\ast} \to Y
\]
from the projectivisation of $R \setminus \{ 0 \}$ to $Y$ which
will be in fact a semigroup homomorphism. The preimage of each
point of $Y$ under the map $\mathrm{\textsf{Init}}$ will give us
an affine space whose codimension $F(\underline{n})$ in
$\mathbb{P}R^{\ast}$ will be determined in Subsect.
\ref{subsection:codimension} in order to show finally an explicit
formula for the generalised Poincar\'e series in terms of the
embedded resolution in Subsect. \ref{subsection:final}.
\medskip

\subsection{Definition of the map $\mathrm{\textsf{Init}}$.} \label{subsection:init}

\begin{Not}
Let be the set $R^{\ast}:=\{z \in R \mid v_i (z) < \infty, ~ 1 \le
i \le r \}$. For $g \in R^{\ast}$, let $\Gamma_g$ be the strict
transform of the curve given by $g$. Since $E_i$ and $\Gamma_g$ have no common
components, the set $E_i \cap \Gamma_g$ is finite for
all $i \in \{1, \ldots
, s \}$. Let $I_0:=\{ \sigma:=(i_1,i_2) \in \{1, 2, \ldots , s
\} \times \{1, 2, \ldots, s \} \mid i_1 < i_2,~ E_{i_1} \cap
E_{i_2} \ne \varnothing \}$. For every $\sigma=(i_1,i_2) \in I_0$,
let $P_{\sigma}:=E_{i_1} \cap E_{i_2}$ and $i_1(\sigma):=i_1,
i_2(\sigma):=i_2$. Let $J_0:=\{1,2,\ldots,r \}$. For $j \in J_0$,
let $P_j:=E_{i_1(j)} \cap \widetilde{C}_j$, where
$\widetilde{C}_j$ is the strict transform of $\pi$ of the
component $C_j$ of the curve $C$, and $E_{i_1(j)}$ is the
component of the exceptional divisor of $\pi$ which intersects
$\widetilde{C}_j$.
\end{Not}
\medskip

We can distinguish three types of points belonging to each
irreducible component $E_i$ of the exceptional divisor of $\pi$,
namely: points of type $P_{\sigma}$ for some $\sigma \in I_0$,
points of type $P_j$ for some $j \in J_0$ and smooth points of the
exceptional divisor.
\medskip

\Nr \label{sec:42} \textbf{Points of type $P_{\sigma}$.} They are defined to be the closed
intersection points $E_{i_1} \cap E_{i_2}$ when
$\sigma=(i_1,i_2)$, for $i_1,i_2 \in \{1,2,\ldots,s \}$ and $i_1 <
i_2$. Take such a point $P_{\sigma}$, for $\sigma \in I_0$. Choose
local coordinates $x_{\sigma},y_{\sigma}$ at $P_{\sigma}$ so that
the components of the total transform of the curve are the
coordinate lines: $E_{i_1(\sigma)}:=\{y_{\sigma}=0 \},
E_{i_2(\sigma)}:=\{x_{\sigma}=0 \}$. We will denote the local ring
at $P_{\sigma}$ by $R_{\sigma}$, and its residue field by
$k_{\sigma}$.
\medskip

The quotient ring $R_{\sigma}/x_{\sigma}R_{\sigma}$ (resp.
$R_{\sigma}/y_{\sigma}R_{\sigma}$) is a discrete valuation ring
with associated valuation $\omega_{x_{\sigma}}$ (resp.
$\omega_{y_{\sigma}}$). We consider the canonical maps
\[
\varphi_{x_{\sigma}}: R_{\sigma} \longrightarrow R_{\sigma} /
x_{\sigma} R_{\sigma}
\]
\[
\varphi_{y_{\sigma}}: R_{\sigma} \longrightarrow R_{\sigma} /
y_{\sigma} R_{\sigma}
\]
Let $g \in R^{\ast}$, and denote by $\gamma_g$ the equation at
$P_{\sigma}$ of the strict transform $\Gamma_g$ of $g$ on $X$
(which must be different from $0$); for $\sigma \in I(g):=
\{\sigma \in I_0 \mid P_{\sigma} \in \Gamma_g \}$, we define
\[
n^{\prime}_{\sigma} (g):= \omega_{x_{\sigma}} \left ( \gamma_g
\mathrm{~mod~} x_{\sigma} R_{\sigma} \right )
\]
\[
n^{\prime \prime}_{\sigma} (g):= \omega_{y_{\sigma}}\left (
\gamma_g \mathrm{~mod~} y_{\sigma} R_{\sigma} \right ).
\]
Note that $n^{\prime}_{\sigma} (g) \ne \infty$ and $n^{\prime
\prime }_{\sigma} (g) \ne \infty$, since $\gamma_g \notin
(x_{\sigma})$ and $\gamma_g \notin (y_{\sigma})$. Note also that
$\varphi_{x_{\sigma}}(y_{\sigma})$ (resp.
$\varphi_{y_{\sigma}}(x_{\sigma})$) is a uniformising parameter of
the discrete valuation ring $R_{\sigma}/x_{\sigma} R_{\sigma}$
(resp. $R_{\sigma}/y_{\sigma} R_{\sigma}$). Therefore we have
\[
\varphi_{x_{\sigma}} (\gamma_g)=\alpha_{\sigma} \cdot
\varphi_{x_{\sigma}} (y_{\sigma})^{n^{\prime}_{\sigma}(g)}
\]
\[
\varphi_{y_{\sigma}} (\gamma_g)=\beta_{\sigma} \cdot
\varphi_{y_{\sigma}} (x_{\sigma})^{n^{\prime \prime}_{\sigma}(g)},
\]
where $\alpha_{\sigma}$ is a unit of
$R_{\sigma}/x_{\sigma}R_{\sigma}$ and $\beta_{\sigma}$ is a unit
of $R_{\sigma}/y_{\sigma} R_{\sigma}$. We can also consider the
image $a_{\sigma}$ of $\alpha_{\sigma}$ and the image $b_{\sigma}$
of $\beta_{\sigma}$ in the residue field $k_{\sigma}$ of
$R_{\sigma}$, which are both different from $0$. The elements
$\alpha_{\sigma}$, $\beta_{\sigma}$, $a_{\sigma}$ and $b_{\sigma}$
depend only on the choice of the equations $x_{\sigma}$ and
$y_{\sigma}$ and also of the equation $\gamma_g$ chosen for the
strict transform of the ideal $gR$. However, it is easy to check
that the ratio $\frac{a_{\sigma}}{b_{\sigma}} \in
k^{\ast}_{\sigma}$ is independent of the generator of the
principal ideal $\gamma_g S_{P_{\sigma}}$. We will denote
$\lambda_{\sigma}(g):=\frac{a_{\sigma}}{b_{\sigma}}$.
\medskip

\Nr \textbf{Points of type $P_j$.} They are the points in the
intersection between $E_i$ and the strict transform of the $j$--th
irreducible component of the curve. Let us take a point $P_j$, for
$j \in J_0$ and choose local coordinates $x_j,y_j$ in a
neighbourhood of $P_j$ in such a way that the components of the
total transform $\pi^{-1}(C)$ of  the curve $C$ are $E_{i_1(j)}=
\{y_j=0 \}$ and $\widetilde{C}_j = \{ x_j=0 \}$, and the set $\{
x_j, y_j \}$ is a regular system of parameters of the local ring
$R_{j}:=\oo_{X,P_j}$. We will use $k_j$ to denote the residue
field of $R_j$. Let $g \in R^{\ast}$. Analogously as done in \ref{sec:42}, for each
element $j$ in $J(g):=\{j \in J_0 \mid P_j \in \Gamma_g
\}$, we define the natural numbers
\[
\widetilde{n}^{\prime}_j (g)=\omega_{x_j} \left ( \gamma_g
\mathrm{~mod~} x_j R_j \right ) \ne \infty
\]
\[
\widetilde{n}^{\prime \prime}_j (g)=\omega_{y_j} \left ( \gamma_g
\mathrm{~mod~} y_j R_j \right ) \ne \infty,
\]
as well as two elements $\alpha_j, \beta_j \in R_j$, whose images
$a_j, b_j$ in the residue field $k_j$ of $R_j$ are different from
$0$ and uniquely determined modulo $x_j$ and $y_j$. We also set $\mu_{j}(g):=\frac{a_{j}}{b_{j}}
\in k^{\ast}_{j}$.
\medskip

\begin{Not}
Let $S$ be a quadratic transform of $R$, let $f \in R$. The strict
transform of $f$ in $S$ will be denoted by $(fR)^S$. This concept
will play an important role in \ref{smooth}, Lemma \ref{lem:surj},
as well as in \ref{nr:56} and Lemma \ref{lem:57}. For a recent
account of the theory we refer the reader to \cite{kiyek}. See also \cite[Sect. 2]{moyano2}.
\end{Not}
\medskip

\Nr \label{smooth} \textbf{Points of type ``smooth''.} They are
the points belonging to $\oei$, for $1 \le i \le s$. Since
the points on $E_i$, $1 \le i \le s$ are closed, if we take $g \in R^{\ast}$ and $\Gamma_g$
the strict transform of the curve given by $g$, then we can count
the number of intersection points between $\Gamma_g$ and $\oei$
with multiplicities as
\[
n_i(g)= \sum_{P \in \oei \cap \Gamma_g} \left ( (gR)^{S_P} \cdot
(\mm R)^{S_P} \right ), ~ ~ ~ \mathrm{~for~} i \in \{1, \ldots ,
s\}.
\]
\medskip

\Nr We want to get a space $Y$ and a map between $\pro R^{\ast}$
and $Y$ so that the preimage of any point of $Y$ is an affine
space of finite codimension, and such a map is a locally trivial
fibration over each connected component of $Y$. This will allow us
to apply certain integration rules to the map in order to connect
it with Definition \ref{poincareseries}.
\medskip

For any (topological) space $\mathfrak{X}$, and $m \in
\mathbb{N}$, we define the $m$--\emph{th symmetric power} of
$\mathfrak{X}$ to be $\mathcal{S}^m \mathfrak{X} := \mathfrak{X}^m
/ S_m$, where $S_m$ is the group of permutations of $m$ elements.
\medskip

\begin{Not} \label{Yn}
Consider the sets $I_0$ and $J_0$. For $I \subset I_0$ and $J
\subset J_0$, let
\begin{eqnarray}
\mathcal{N}_{I,J}& := &  \{ \underline{n}:= (n_i,
n_{\sigma}^{\prime}, n_{\sigma}^{\prime \prime},
\widetilde{n}_{j}^{\prime}, \widetilde{n}_{j}^{\prime \prime})
\mid n_i \ge 0,~1 \le i \le s;  ~ n^{\prime}_{\sigma}>0,  \nonumber \\
&  &~ n^{\prime \prime}_{\sigma}>0,~ \sigma \in
I;~\widetilde{n}^{\prime}_{j}>0,~\widetilde{n}^{\prime
\prime}_{j}>0, ~ j \in J  \}. \nonumber
\end{eqnarray}
\medskip

For each $\underline{n} \in \mathcal{N}_{I,J}$, we define
\[
Y_{\underline{n}} := \prod_{i=1}^{s} \mathcal{S}^{n_i}
\overset{\circ}{E_i} \times \prod_{\sigma \in I}
{k}_{\sigma}^{\ast} \times \prod_{j \in J} {k}_{j}^{\ast},
\]
with $\mathcal{S}^{n_i} \overset{\circ}{E_i}$ being the $n_i$--th
symmetric power of $\overset{\circ}{E_i}$, $k^{\ast}_{\sigma}$ the
residue field at $P_{\sigma}$ minus $0$ corresponding to the
element $\sigma$ of $I$, and $k^{\ast}_j$ the residue field minus
zero at $P_j$ corresponding to the element $j$ of $J$. Setting
\[
Y:= \bigcup_{\substack{I \subset I_0 \\ J \subset J_0}}
\bigcup_{\underline{n} \in \mathcal{N}_{I,J}} Y_{\underline{n}},
\eqno(\dag)
\]
as in \cite[p.~203]{cadegu11}, it is easily checked that $Y$  can
be endowed with a structure of semigroup.
\end{Not}
\medskip

\Nr Recall that $R^{\ast}= \{z \in R \mid v_i (z) < \infty, ~1 \le i
\le r \}$. Consider the quotient $R^{\ast}/\sim$, where ``$\sim$"
denotes the following equivalence relation: For two elements  $a,
b \in R^{\ast}$, we say that $a \sim b$ if there exists an element
$u \in k_R \cong K$ such that $a=ub$. The set $
\mathbb{P} R^{\ast} := R^{\ast} / \sim $ can be endowed with a structure
of semigroup just considering the multiplication of functions.
Define a semigroup homomorphism
\begin{displaymath}
\begin{array}{rccc}
\mathrm{\textsf{Init}}: & \mathbb{P}R^{\ast} & \longrightarrow & Y
= \bigcup_{\substack{I \subset I_0 \\ J \subset J_0}}
\bigcup_{\underline{n} \in \mathcal{N}_{I,J}} Y_{\underline{n}}
\end{array}
\end{displaymath}
between $\mathbb{P} R^{\ast}$ and $Y$, in which the image
$\mathrm{\textsf{Init}}(g)$ of any element $g \in \mathbb{P}
R^{\ast}$ is an element of the connected component
$Y_{\underline{n}}$ of $Y$ corresponding to $I(g),J(g)$ and with
$\underline{n}=\left (n_i(g), n^{\prime}_{\sigma} (g),
n_{\sigma}^{\prime \prime} (g), \widetilde{n}_j^{\prime} (g),
\widetilde{n}_j^{\prime \prime} (g) \right )$, which is defined in
every factor (connected component) of $Y_{\underline{n}}$ as
follows:
\begin{enumerate}
    \item[--] in $\mathcal{S}^{n_i} \oei$, $1 \le i \le s$, $\mathrm{\textsf{Init}}(g)$
    is represented by the set of intersection points between
    $\Gamma_g$ and $\oei$  counted with multiplicities;
    \item[--] in $k_{\sigma}^{\ast}$, $\sigma \in I(g)$,
    $\mathrm{\textsf{Init}}(g)$ is represented by the quotient
    $\frac{a_{\sigma}(g)}{b_{\sigma}(g)}$;
    \item[--] in $k_{j}^{\ast}$, $j \in J(g)$,
    $\mathrm{\textsf{Init}}(g)$ is represented by the quotient
    $\frac{a_{j}(g)}{b_{j}(g)}$.
\end{enumerate}
\medskip

We take the opportunity to prove the following key result (it was
already omitted for the complex case in \cite{cadegu11}).
\medskip

\begin{Lem} \label{lem:surj}
The map \textsf{Init} is surjective.
\end{Lem}

\dem~ Let us fix $I \subset I_0$, $J \subset J_0$, $\underline{n}
\in \mathcal{N}_{I,J}$. Consider an element
\[
y_{\underline{n}} = \left ( (\{P_{i,j} \}_{j=1}^{n_i};~i=1, \ldots
, s), (\lambda_{\sigma};\sigma \in I), (\mu_j; j \in J) \right )
\in Y_{\underline{n}}.
\]
We will construct $g \in R^{\ast}$ such that
$\mathrm{\textsf{Init}}(g)=y_{\underline{n}}$ step by step, in
fact it suffices to show the surjectivity in each factor of the
space $Y_{\underline{n}}$.
\medskip

For a point $P \in \oei$, Lemma 3.1.11 of \cite{moyano} (see also
\cite[Sect. 6]{moyano2}) ensures the existence of an irreducible
element $g_P \in R$ such that the strict transform of $g_P R$ in
$X$ is smooth and transversal to the exceptional divisor $E$ at
the point $P \in \oei$, i.e., $(g_P R)^{S_Q}=S_Q$ for all $Q \ne
P$ ($S_Q$ denotes the local ring of $X$ at the point $Q \in E$)
and $(g_P R)^{S_P}=\gamma_P S_P$ so that $\{x_P,\gamma_P \}$ is a
parameter system of $S_P$ (here $x_P$ is a local equation for
$E_i$ at $P$; i.e., $x_P S_P = (\mm_R R)^{S_P}$). If we take
\[
\left ( \{P_{ij} \}_{j=1}^{n_i};~i=1, \ldots , s \right ) \in
\prod_{i=1}^s \mathcal{S}^{n_i} \oei,
\]
then it is easily seen that $g^{(1)}=\prod_{i,j} g_{P_{ij}} \in R$
defines an element of $R$ such that $I(g^{(1)})=J(g^{(1)}) =
\varnothing$, and moreover $n_i (g^{(1)})=n_i$ for $1 \le i \le s$.
\medskip

Let be now $\sigma \in I$, $n_{\sigma}^{\prime},
n_{\sigma}^{\prime \prime} >0$ and $\lambda_{\sigma} \in
k^{\ast}_{\sigma}$. Let us consider the discrete valuation
$\omega_{x_{\sigma}}$ (resp. $\omega_{y_{\sigma}}$) associated
with the ring $R_{\sigma}/ x_{\sigma} R_{\sigma}$ (resp.
$R_{\sigma}/ y_{\sigma} R_{\sigma}$) (Recall that $\{x_{\sigma},
y_{\sigma} \}$ are the equations of the components of the
exceptional divisor at the point $P_{\sigma}=E_{i_1(\sigma)} \cap E_{i_2
(\sigma)}$). Let $u_{\sigma} \in R_{\sigma}$ be a unit such that
$u_{\sigma} + \mm_{\sigma} = \lambda_{\sigma} \in
k_{\sigma}^{\ast}$ and let us take $\gamma_1^{\sigma}=x_{\sigma} +
u_{\sigma} \cdot y_{\sigma}^{n^{\prime}_{\sigma}-1}$,
$\gamma_2^{\sigma}=y_{\sigma} + x_{\sigma}^{n^{\prime
\prime}_{\sigma}-1} \in R_{\sigma}$. Since $\{y_{\sigma},
\gamma_{1}^{\sigma} \}$ is a regular system of parameters, we have
that $\gamma_1^{\sigma} R_{\sigma} \cap R$ is a prime ideal of
height $1$, therefore it is principal and generated by
$g_1^{\sigma} \in R$ whose strict transform in $R_{\sigma}$ is
just $\gamma_1^{\sigma} R_{\sigma}$ (see \cite[Prop. 6.2]{moyano2}). Identically we can argue with
$\gamma_2^{\sigma}$ to get an irreducible element
$g_2^{\sigma} \in R$ with $(g_2^{\sigma} R
)^{R_{\sigma}}=\gamma_2^{\sigma} R_{\sigma}$. Let us denote
$g^{\sigma}=g_1^{\sigma} \cdot g_2^{\sigma}$ and
$\gamma^{\sigma}=\gamma_1^{\sigma} \cdot \gamma_2^{\sigma}$.
Notice that $\gamma^{\sigma}$ is the strict transform of
$g^{\sigma}$. Then we have
\[
\varphi_{x_{\sigma}}(\gamma^{\sigma})= \varphi_{x_{\sigma}} (
u_{\sigma}) \cdot
\varphi_{x_{\sigma}}(y_{\sigma})^{n^{\prime}_{\sigma}}
\]
\[
\varphi_{y_{\sigma}}(\gamma^{\sigma})=
\varphi_{y_{\sigma}}(x_{\sigma})^{n^{\prime \prime}_{\sigma}},
\]
and since $\varphi_{x_{\sigma}}(u_{\sigma})$ is a unit, we have
$n^{\prime}_{\sigma} (g^{\sigma})=n^{\prime}_{\sigma}$ and
$n^{\prime \prime}_{\sigma} (g^{\sigma})=n^{\prime
\prime}_{\sigma}$. Moreover, $a_{\sigma}=\lambda_{\sigma}$,
$b_{\sigma}=1$ and therefore
$\frac{a_{\sigma}}{b_{\sigma}}=\lambda_{\sigma} \in
k_{\sigma}^{\ast}$.
\medskip

Thus let $g^{(2)}=\prod_{\sigma \in I} g^{\sigma} \in R$. From
the previous construction, it follows that $n_i (g^{(2)})=0$ for
all $i=1, \ldots , s$; $J(g^{(2)})= \varnothing$; $I(g^{(2)})=I$ and
moreover, for every $\sigma \in I$, we have
$n^{\prime}_{\sigma}(g^{(2)})=n^{\prime}_{\sigma}$, $n^{\prime
\prime}_{\sigma} (g^{(2)})^=n^{\prime \prime}_{\sigma}$ and
$\lambda_{\sigma} (g^{(2)})=\lambda_{\sigma} \in
k_{\sigma}^{\ast}$.
\medskip

For every $j \in J$ (and consequently for $J$) we can proceed
analogously to the previous case in order to get an element
$g^{(3)} \in R$ such that $n_i (g^{(3)})=0$ for all $i= 1, \ldots
, s$; $I(g^{(3)})=\varnothing$; $J(g^{(3)})=J$ and
$\widetilde{n}^{\prime}_{j} (g^{(3)})=\widetilde{n}^{\prime}_{j}$,
$\widetilde{n}^{\prime \prime}_{j} (g^{(3)})=\widetilde{n}^{\prime
\prime}_{j}$ and $\mu_j (g^{(3)})=\mu_j \in k_{j}^{\ast}$.
\medskip
\vspace{10cm}

As a consequence, we obtain an element $g=g^{(1)} \cdot g^{(2)}
\cdot g^{(3)}$ satisfying
\[
\mathrm{\textsf{Init}}(g)= \left ( (\{P_{i,j} \}_{j=1}^{n_i};~i=1,
\ldots , s), (\lambda_{\sigma};\sigma \in I), (\mu_j; j \in J)
\right ) \in Y_{\underline{n}}.
\]
\qed
\medskip

\begin{Not}
Take into account that $X$ is an irreducible integral scheme, and $E_i$ is an
irreducible closed subset of $X$ for every $1 \le i \le s$.
Consider $(E_i, \oo_{E_i}) \subset (X, \oo_X)$ with reduced
structural sheaf. Let $\eta_i$ be the generic point of $E_i$. We consider the local ring
$\oo_{X,E_i}$ of $X$ along $E_i$, which is nothing but the ring
$\oo_{X,\eta_i}$. Then there exists an affine open subset $U$ in
$X$ with $\eta_i \in U$. Let be the ring $A:=\Gamma (U, \oo_X)$.
Then the generic point $\eta_i$ corresponds to an ideal
$\mathfrak{p}_i \in \mathrm{Spec}(A)$, and therefore
$\oo_{X,E_i}=A_{\mathfrak{p}_i}$. It is easy to check that
$A_{\mathfrak{p}_i}$ is a discrete valuation ring, with associated
discrete valuation $w_i$, for $1 \le i \le s$.
\end{Not}
\medskip

\begin{Prop} \label{lemma:clave}
For $g, g^{\prime} \in R^{\ast}$, then
$\mathrm{\textsf{Init}}(g)=\mathrm{\textsf{Init}}(g^{\prime})$ if
and only if $g = \zeta g^{\prime} + h$, where $\zeta \in k_R
\setminus \{ 0 \}$, $w_i (h)
> w_i (g)=w_i (g^{\prime})$ for all $1 \le i \le s$ and
$v_{j} (h) > v_{j} (g)=v_{j} (g^{\prime})$ for all $j \in \{1,
\ldots , r \}$.
\end{Prop}

\dem~Let us take a point $q \in E_i \cap U \ne \varnothing$. The
point $q$ corresponds to a prime ideal $\mathfrak{q} \in
\mathrm{Spec}(A)$ with $\mathfrak{q} \supset \mathfrak{p}_i$.
Therefore we have $\oo_{X,q}=A_{\mathfrak{q}}$ and
$\oo_{X,E_i}=A_{\mathfrak{p}_i}= \left (A_{\mathfrak{q}} \right
)_{\mathfrak{p}_i A_{\mathfrak{q}}}$.
Note that the field of rational functions of $X$ (resp.~of
$E_i$) is $K(X):=\mathrm{Quot}(\oo_{X,E_i})=\mathrm{Quot}(R)$ (resp.
$K(E_i):= A_{\mathfrak{p}_i}/\mathfrak{p}_i A_{\mathfrak{p}_i}$).
\medskip

We consider $g,g^{\prime} \in R$, the morphism $\rho: X \to U$ and
the induced ring homomorphism $\theta_i: \Gamma (U,\oo_X) \to
\Gamma (U \cap E_i, \oo_{E_i})$, for every $1 \le i \le s$. Then
$\rho (g) \in \Gamma (U, \oo_X)$ and $\theta_i (\rho (g)) \in
\Gamma (U \cap E_i, \oo_{E_i})$. Denote by $[\theta_i (\rho (g))]$
the class of $\theta_i (\rho (g))$ as a rational function on
$E_i$. Let $\iota_i: E_i \hookrightarrow X$ be the inclusion map.
Define $\widetilde{g}:=g \circ \pi$ and
$\widetilde{g^{\prime}}:=g^{\prime} \circ \pi$ in $\Gamma (X,
\oo_{X})$. Then $\widetilde{g} \circ \iota_i =[\theta_i
(\rho(g))]$, $\widetilde{g^{\prime}} \circ \iota_i =[\theta_i
(\rho(g^{\prime}))]$ and we can associate with $\widetilde{g}
\circ \iota_i$ (resp. $\widetilde{g^{\prime}} \circ \iota_i$) a
divisor
\[
\mathrm{div} (\widetilde{g} \circ \iota_i)=\sum_{P \in E_i} \left
( \Gamma_g \cdot E_i \right )_P \cdot P
\]
(resp. $\mathrm{div}(\widetilde{g^{\prime}} \circ \iota_i)=\sum_{P
\in E_i} \left ( \Gamma_{g^{\prime}} \cdot E_i \right )_P \cdot
P$), where $\Gamma_g$ and $\Gamma_{g^{\prime}}$ are the strict
transforms of $g$ and $g^{\prime}$ by $\pi$, respectively.
\medskip

Set $\psi:=\frac{\widetilde{g}}{\widetilde{g^{\prime}}} \in
\oo_{X, \eta_i}$, for every $1 \le i \le s$. The condition
$\mathrm{\textsf{Init}}(g)=\mathrm{\textsf{Init}}(g^{\prime})$
means, in particular, that $\mathrm{div} (\widetilde{g} \circ
\iota_i)= \mathrm{div} (\widetilde{g^{\prime}} \circ \iota_i)$,
and then the restriction of $\psi$ to $E_i$, i.e., the quotient
$\frac{[\theta_i(\rho(g))]}{[\theta_i(\rho(g^{\prime}))]}$, is a
regular function on $E_i$. Similarly, since the irreducible
components $\widetilde{C}_j$ of the strict transform of the curve
$C$ are smooth, we can also show that the restriction of $\psi$ to
each $\widetilde{C}_j$ is a regular function on $\widetilde{C}_j$,
for every $1 \le j \le r$.
\medskip

Moreover, regular functions defined over projective lines are constant, so the 
function $\psi$ is constant on each component $E_i$ of the
exceptional divisor. In fact, as $E_i=\mathrm{Proj} \left ( k_i
[\overline{x},\overline{y}] \right )$ for some finite field
extension $k_i \supseteq k_R$ (according to Sect.~\ref{sec:2}),
we have
\[
\frac{[\theta_i(\rho(g))]}{[\theta_i (\rho (g^{\prime}))]}=\zeta_i
\in k_i \setminus \{ 0 \}, ~ ~ ~ ~ ~ ~ ~ 1 \le i \le s.
\]
We see now that, for every $1 \le i \le s$, there
exists a constant $\zeta \in k_R \setminus \{ 0 \}$ such that
$\zeta_i=\zeta$ for all $1 \le i \le s$.
\medskip

Let $P_{\sigma}$, $\sigma=(i_1,i_2)$ be the intersection point
between $E_{i_1}$ and $E_{i_2}$, for $i_1 < i_2$. Let
$R_{\sigma}:=\oo_{X, P_{\sigma}}$ be the local ring of $X$ at
$P_{\sigma}$, and consider the local equation $\gamma_g$ (resp.
$\gamma_{g^{\prime}}$) of $\Gamma_g$ (resp. $\Gamma_{g^{\prime}}$)
at $P_{\sigma}$, and the local equations $\{y_{\sigma}=0 \}$ and
$\{x_{\sigma}=0 \}$ of $E_{i_1}$ and $E_{i_2}$, respectively. Let
\begin{displaymath}
\begin{array}{lccc}
\varphi_{x_{\sigma}}: & R_{\sigma} & \longrightarrow & R_{\sigma}/x_{\sigma} R_{\sigma}\\
& \gamma_g &
\mapsto & \gamma_g \mathrm{~mod~} x_{\sigma},\\
\end{array}
\end{displaymath}
\begin{displaymath}
\begin{array}{lccc}
\varphi_{y_{\sigma}}: & R_{\sigma} & \longrightarrow & R_{\sigma}/y_{\sigma} R_{\sigma}\\
& \gamma_g &
\mapsto & \gamma_g \mathrm{~mod~} y_{\sigma}.\\
\end{array}
\end{displaymath}
be the canonical homomorphisms. From the assumption
$\mathrm{\textsf{Init}}(g)=\mathrm{\textsf{Init}}(g^{\prime})$, it
follows that
\[
\deg_X \left ( \Gamma_g \cdot E_i \right ) = \deg_X \left (
\Gamma_{g^{\prime}} \cdot E_i \right ) = n_{\sigma}
\]
\[
\deg_X \left ( \Gamma_g \cdot E_j \right ) = \deg_X \left (
\Gamma_{g^{\prime}} \cdot E_j \right ) = m_{\sigma}.
\]
Then there exist $\lambda_{\sigma}, \lambda^{\prime}_{\sigma},
\mu_{\sigma}, \mu^{\prime}_{\sigma}$ units of $R_{\sigma}$, so
that
\[
\gamma_g - \lambda_{\sigma} y_{\sigma}^{n_{\sigma}},
\gamma_{g^{\prime}}-\lambda^{\prime}_{\sigma}
y_{\sigma}^{n_{\sigma}} \in x_{\sigma} R_{\sigma}
\]
\[
\gamma_{g}-\mu_{\sigma} x_{\sigma}^{m_{\sigma}},
\gamma_{g^{\prime}}-\mu^{\prime}_{\sigma} x_{\sigma}^{m_{\sigma}}
\in y_{\sigma} R_{\sigma}.
\]
Therefore
\[
\varphi_{x_{\sigma}}(\gamma_{g})= \varphi_{x_{\sigma}}
(\lambda_{\sigma}) \cdot
\varphi_{x_{\sigma}}(y_{\sigma})^{n_{\sigma}}
\]
\[
\varphi_{x_{\sigma}}(\gamma_{g^{\prime}})= \varphi_{x_{\sigma}}
(\lambda^{\prime}_{\sigma}) \cdot
\varphi_{x_{\sigma}}(y_{\sigma})^{n_{\sigma}}
\]
\[
\varphi_{y_{\sigma}}(\gamma_{g})= \varphi_{y_{\sigma}}
(\mu_{\sigma}) \cdot \varphi_{y_{\sigma}}(x_{\sigma})^{m_{\sigma}}
\]
\[ \varphi_{y_{\sigma}}(\gamma_{g^{\prime}})= \varphi_{y_{\sigma}}
(\mu^{\prime}_{\sigma}) \cdot
\varphi_{y_{\sigma}}(x_{\sigma})^{m_{\sigma}}.
\]
Thus the value of $\frac{[\theta_{i_1} (\rho (g))]}{[\theta_{i_1}
(\rho (g^{\prime}))]}$ at the point $P_{\sigma}$ is just
$\frac{\varphi_{y_{\sigma}}(\mu_{\sigma})}{\varphi_{y_{\sigma}}(\mu^{\prime}_{\sigma})}=
\frac{b_{\sigma}}{b^{\prime}_{\sigma}} = \zeta_{i_1}$, and the
value of $\frac{[\theta_{i_2} (\rho (g))]}{[\theta_{i_2} (\rho
(g^{\prime}))]}$ at $P_{\sigma}$ is
$\frac{\varphi_{x_{\sigma}}(\lambda_{\sigma})}{\varphi_{x_{\sigma}}(\lambda^{\prime}_{\sigma})}
= \frac{a_{\sigma}}{a^{\prime}_{\sigma}} =\zeta_{i_2}$. Since
$\mathrm{\textsf{Init}}(g) = \mathrm{\textsf{Init}} (g^{\prime})$,
we get
\[
\frac{\varphi_{x_{\sigma}}(\lambda_{\sigma})}{\varphi_{y_{\sigma}}(\mu_{\sigma})}
=
\frac{\varphi_{x_{\sigma}}(\lambda^{\prime}_{\sigma})}{\varphi_{y_{\sigma}}(\mu^{\prime}_{\sigma})}
\Longleftrightarrow
\frac{\varphi_{y_{\sigma}}(\mu_{\sigma})}{\varphi_{y_{\sigma}}(\mu^{\prime}_{\sigma})}
=
\frac{\varphi_{x_{\sigma}}(\lambda_{\sigma})}{\varphi_{x_{\sigma}}(\lambda^{\prime}_{\sigma})}
\Longleftrightarrow \zeta_{i_1}=\zeta_{i_2}.
\]
Furthermore, taking into account that $E_1 \cup \ldots \cup E_s$
is a connected set, we obtain that $\psi$ is a regular function on
the total transform $\pi^{-1}(C)$, and it is equal to a constant
$\zeta \ne 0$ on the exceptional divisor, which must be an element
belonging to $k_R$, because there is (at least) one $E_i$ which is
a scheme $\mathrm{Proj}(k_R[\overline{x},\overline{y}])$.
\medskip

Now we can take the function $h:=g-\zeta g^{\prime}$. For each $1
\le i \le s$, we consider the canonical homomorphism
\begin{displaymath}
\begin{array}{lccc}
\phi_i: & A_{\mathfrak{p}_i} & \longrightarrow & A_{\mathfrak{p}_i}/ \mathfrak{p}_i A_{\mathfrak{p}_i}=K(E_i)=k_i \\
& a &
\mapsto & a \mathrm{~mod~} \mathfrak{p}_i.\\
\end{array}
\end{displaymath}
The quotient $\frac{g}{g^{\prime}} \in \oo_{X,
\eta_i}=\oo_{X,E_i}=A_{\mathfrak{p}_i}$, and therefore $\phi_i
\left ( \frac{g}{g^{\prime}} \right ) = \zeta \neq 0$, then
$\frac{g}{g^{\prime}} - \zeta \in
\mathrm{ker}(\phi_i)=\mathfrak{p}_i A_{\mathfrak{p}_i}$ and
\[
0 < w_i \left ( \frac{g}{g^{\prime}} - \zeta \right ) = w_i \left
( \frac{g- \zeta g^{\prime}}{g^{\prime}} \right ) = w_i (g - \zeta
g^{\prime}) - w_i (g^{\prime}).
\]
Hence $w_i(h)>w_i(g)=w_i(g^{\prime})$ for all $1 \le i \le s$. We can repeat the same
argument for
the valuations $v_j$, $1 \le j \le r$; the intersection $C_j \cap \mathrm{Spec}(A)$ is not
empty, and so there exists an ideal $\mathfrak{q}_j \in
\mathrm{Spec}(A)$ which is the generic point of $C_j$ and
$\oo_{X,C_j}=A_{\mathfrak{q}_j}$, which is a discrete valuation
ring with associated discrete valuation $v_j$, for $1 \le j \le
r$. Then it is enough to take the homomorphism
\[
\phi_j: A_{\mathfrak{q}_j} \longrightarrow A_{\mathfrak{q}_j} /
\mm (A_{\mathfrak{q}_j}),
\]
where $\mm (A_{\mathfrak{q}_j})$ is the maximal ideal of the local
ring $A_{\mathfrak{q}_j}$.
\medskip

Conversely, if $g=\zeta g^{\prime}+h$ with $\zeta \in k_R
\setminus \{0\}$, $w_i (g - \zeta g^{\prime}) > w_i (g) =
w_i(g^{\prime})$ for $1 \le i \le s$ and $v_j (g-\zeta g^{\prime})
> v_j(g) = v_j (g^{\prime})$ for all $1 \le j \le r$, then the
function $\psi$ is regular on every component $E_i$ of $E$ for $1
\le i \le s$ and on every component $\widetilde{C}_j$ of the
strict transform of the curve, and therefore $\psi |_{E} \equiv
\zeta$. Hence $\mathrm{div}(\widetilde{g^{\prime}} \circ \iota_i)
= \mathrm{div}(g \circ \iota_i)$ for $1 \le i \le s$ and the
intersection points of the strict transforms $\Gamma_g$, and
$\Gamma_{g^{\prime}}$ with each component $E_i$ of $E$ coincide
(counting multiplicities). For $1 \le j \le r$, we also have
\[
\left ( \Gamma_{g} \cdot \widetilde{C}_j \right )_{P_j} = v_j (g)
- w_{i_1(j)} (g) = v_j (g^{\prime}) - w_{i_1(j)} (g^{\prime}) =
\left ( \Gamma_{g^{\prime}} \cdot \widetilde{C}_j \right )_{P_j}.
\]
Moreover, for points of type $P_{\sigma}$ (the same holds for
points of type $P_j$) and following previous notations, the value
of $\frac{[\theta_i (\rho (g))]}{[\theta_i (\rho (g^{\prime}))]}$
at $P_{\sigma}$ is $\zeta =
\frac{\varphi_{x_{\sigma}}(\lambda_{\sigma})}{\varphi_{x_{\sigma}}(\lambda^{\prime}_{\sigma})}=
\frac{a_{\sigma}}{a^{\prime}_{\sigma}}$, but it is also $\zeta =
\frac{b_{\sigma}}{b^{\prime}_{\sigma}}$, then
$\frac{a_{\sigma}}{a^{\prime}_{\sigma}}=\frac{b_{\sigma}}{b^{\prime}_{\sigma}}$.
\qed
\medskip

\begin{Cor} \label{cor:F(n)}
For any point $\varsigma \in Y_{\underline{n}} \subset Y$, the
preimage $\mathrm{\textsf{Init}}^{-1}(\varsigma)$ consists of an
affine space given by a point $g$ plus the ideal
$I_{\underline{n}}$, where
\[
I_{\underline{n}} = \left \{z \in R^{\ast} \mid w_i (z) > w_i
(\underline{n}), ~ 1 \le i \le s; ~v_j (z) > v_j (\underline{n}),~
1 \le j \le r \right \}.
\]
This space has finite codimension $F(\underline{n})$ in $\pro R^{\ast}$.
Furthermore, over each connected component $Y_{\underline{n}}$ of
$Y$, the map $\mathrm{\textsf{Init}}$ is a locally trivial
fibration.
\end{Cor}
\medskip

\subsection{Computation of the codimension F(\underline{n}).}
\label{subsection:codimension}

\Nr \label{nr:413} Let $g \in R^{\ast}$ and $w_i$ the divisorial
valuation corresponding to the component $E_i$ of $\pi$, for all
$1 \le i \le s$. Set $\underline{w}(g):=(w_1 (g), \ldots ,
w_s(g))$. We know that the value $w_i (g)$ is the multiplicity of
the lifting $g \circ \pi$ of $g$ along $E_i$, for all $1 \le i \le
s$. Let $\mathcal{D}(g)$ be the total transform by $\pi$ of the
curve defined by the ideal $gR$; i.e. $\mathcal{D}(g)$ is the
divisor in $X$ given by
\[
\mathcal{D}(g) = \sum_{i=1}^s w_i(g) E_i + \Gamma_g,
\]
where $\Gamma_g$ is the strict transform of the curve defined by
$gR$. We have now
\[
0 = \deg_{X}(\mathcal{D}(g)\cdot E_i) = \deg_{X}(\Gamma_g\cdot
E_i) + \sum_{j=1}^s w_j(g) \deg_{X}(E_i\cdot E_j)
\]
and so $\widehat{n}_i(g) :=  \deg_{X}(\Gamma_g\cdot E_i) = -
\sum_{j=1}^s w_j(g)n_{ij}$, being $N=(n_{ij})$ the intersection
matrix (see \ref{nr:matrices}). Thus we have that
\[
\underline{w} (g) := (w_1(g), \ldots, w_s(g)) = (\widehat{n}_1(g),
\ldots, \widehat{n}_s(g)) \cdot M,
\]
where $M:=-N^{-1}$. The above formula provides the precise
relation between the values $\underline{w} (g) := (w_1(g), \ldots, w_s(g))$
and $\widehat{\underline{n}}(g):=(\widehat{n}_1(g), \ldots,
\widehat{n}_s(g))$.
\medskip

Let us now fix $\widehat{\underline{n}} = (\widehat{n}_1,\ldots,
\widehat{n}_s)\in \mathbb{Z}_{\ge 0}^s$ and the corresponding
divisor $A= \sum_{i=1}^s w_i E_i$, where $\underline{w} =
(w_1,\ldots, w_s) = \widehat{\underline{n}} \cdot M$. Let
$J^{D}(\underline{w})$ be the divisorial ideal of $R$ defined by
\[
J^{D}(\underline{w}) = \left \{z \in R^{\ast} \mid w_i (z) \ge
w_i,
 ~ 1 \le i \le s \right \}
\]
and $h^{D}(\underline{w}) = \dim_{k_R} R/J^{D}(\underline{w})$.
\medskip

The next proposition is known as the Hoskin--Deligne formula:
\medskip

\begin{Prop}\label{H-De}
Let $\mathbb{K} = \sum_{i=1}^{s} E_i^{\ast}$ be the canonical
divisor on $X$ and consider $A= \sum_{i=1}^s \alpha_i E^{*}_i$ be the
expression of the divisor $A$ on the basis $\{E^*_i\}$. Then
\[
h^{D}(\underline{w}) =  \frac 12 \sum_{i=1}^s h_i \alpha_i
(\alpha_i+1)\; .
\]
\end{Prop}

\dem~ The genus formula gives us that
\[
h^0(\oo_{A}) = h^{D}(\underline{w}) = -  \frac{\deg_X (A \cdot A)
+ \deg_X (A \cdot \mathbb{K})}{2}\; .
\]
Now, by Equation (\ddag) in \ref{nr:215} we deduce
\[
\deg_X (A \cdot A) = - \sum_{i=1}^{s} \alpha_i^2 h_i
\]
\[
\deg_X (A \cdot \mathbb{K}) = - \sum_{i=1}^{s} \alpha_i h_i.
\]
Then the codimension of the ideal $J^{D}(\underline{w})$ is equal
to $\frac{1}{2} \sum_{i=1}^{s} h_i \alpha_i (\alpha_i +1)$. \qed
\medskip

The point is now to compute the intersection degrees $\deg_X (A
\cdot A)$ and $\deg_X (A \cdot \mathbb{K})$ in terms of the matrix $M$ and
$\widehat{n}_1, \ldots, \widehat{n}_s$.
\medskip

Let $\overset{\bullet}{\nu_i}$ (resp. $\overset{\circ}{\nu_i}$) be
the number of components $E_j$ of $E$ intersecting $E_i$ for $j
\ne i$ (resp. this number plus the number of strict transforms of
the curve $C$ intersecting $E_i$), but counting this number of
components as many times as the degree of the extension
$h_i=[k_i:k_R]$ says. Notice also that $\overset{\bullet}{\nu_i} =
\deg_X(E_i \cdot (\sum_{j\neq i}E_j)) = \sum_{j\neq i} n_{ij}$.
\medskip

\begin{Lem} \label{lem:DK}
For $i=1, \ldots, s$, let $\varepsilon_i:=2 h_i -
\overset{\bullet}{\nu_i}$. Then we have
\[
\deg_X (A \cdot \mathbb{K}) = \widehat{\underline{n}} \cdot
\underline{1}^t - \widehat{\underline{n}} \cdot M \cdot
\underline{\varepsilon}^t.
\]
\end{Lem}

\dem~ For every $1 \le i \le s$, we claim that
\[
 \deg_X (E_i \cdot E_i) = -2 h_i -\deg_X (\mathbb{K} \cdot E_i).
 \eqno(*)
\]

In fact, let $P=(p_{ij})_{1 \le i,j \le s}$ be the proximity
matrix of $\pi$ (cf. \ref{nr:matrices}). From the definitions of
intersection matrix (see again \ref{nr:matrices}) and the divisor
$\mathbb{K}$, we have that $\deg_X (E_i \cdot E_i) =
-\sum_{j=1}^{s} p_{ij}^2 h_i$ and also $\deg_X (\mathbb{K} \cdot E_i) =
- \sum_{j=1}^{s} p_{ij} h_i$, therefore $\deg_X (E_i \cdot E_i) +
\deg_X (\mathbb{K} \cdot E_i) = - \sum_{j=1}^{s} h_i p_{ij}
(p_{ij} +1)$. For $1 \le j \le s$, $p_{ij}(p_{ij}+1)=0$ if and
only if $p_{ij}=0$ or $p_{ij}=-1$, which occurs whenever $i \ne
j$. Then, for every  $1 \le i \le s$ we have
$$
\deg_X (E_i \cdot E_i) + \deg_X (\mathbb{K} \cdot E_i) = -h_i
(p_{ii}^2 + p_{ii}) = -h_i (1+1) = - 2 h_i
$$
and Equation ($\ast$) has been shown. From this follows
\begin{equation}
\sum_{j=1}^{s} p_{ij} h_i = - \deg_X (\mathbb{K} \cdot E_i) = 2h_i
+ \deg_X (E_i \cdot E_i), ~ ~
 1 \le i \le s. \nonumber
\end{equation}

On the other hand, the definition of the intersection matrix
implies
\begin{equation} \nonumber
N \cdot \begin{pmatrix}
  1 \\
  \vdots \\
  1
\end{pmatrix}
=
\begin{pmatrix}
  \deg_X (E_1 \cdot E_1) + \overset{\bullet}{\nu_1} \\
  \vdots \\
  \deg_X (E_s \cdot E_s) + \overset{\bullet}{\nu_s}
\end{pmatrix}
\end{equation}

If we write $\underline{1} := (1, 1, \ldots, 1)$ and
$\underline{\varepsilon}:= (\varepsilon_1, \ldots ,
\varepsilon_s)$, with $\varepsilon_i:= 2 h_i -
\overset{\bullet}{\nu_i} $, for $1 \le i \le s$, then we obtain
\[
N \cdot \underline{1}^t + \underline{\varepsilon}^t = P \cdot
\Delta \cdot \underline{1}^t,
\]
where $\Delta$ is the diagonal matrix defined in
\ref{nr:matrices}. Therefore
\begin{align*}
 - \deg_X (A \cdot \mathbb{K}) & = (\alpha_1, \ldots ,
\alpha_s) \cdot \Delta \cdot (1, \ldots , 1)^t  \\
  & = \underline{w} \cdot P \cdot \Delta \cdot
(1, \ldots , 1)^t \\
 & = \widehat{\underline{n}} \cdot M \cdot P \cdot \Delta \cdot
\underline{1}^t \\
 & = \widehat{\underline{n}} \cdot M \cdot \left ( N \cdot
\underline{1}^t + \underline{\varepsilon}^t \right ) \\
  & = -\widehat{\underline{n}} \cdot \underline{1}^t +
\widehat{\underline{n}} \cdot M \cdot \underline{\varepsilon}^t.
\end{align*}
\qed
\medskip

\begin{Rem}
Let $\beta_i$ be the number of components in $\{E_j: j\neq i\}$
such that $E_i\cap E_j\neq \varnothing$. Then
$\overset{\bullet}{\nu_i} = h_i \beta_i$ and so $\varepsilon_i =
h_i (2- \beta_i)$. In the complex case $2-\beta_i$ is the Euler
characteristic of the space $\overset{\bullet}{E_i} = E_i
\setminus \bigcup_{j\neq i} E_j$; thus in some sense
$\varepsilon_i$ may be interpreted as the ``Euler characteristic"
of $\overset{\bullet}{E_i}$.
\end{Rem}
\medskip

\begin{Lem}\label{lem:DD}
We have
\[
- \deg_X (A \cdot A)= \widehat{\underline{n}} \cdot M \cdot
\widehat{\underline{n}}^t.
\]
\end{Lem}

\dem~Since $M=-N^{-1}$, $\underline{w}
=\widehat{\underline{n}} \cdot M$, and $M$ is a symmetric
matrix, we obtain
\begin{align*}
 - \deg_X (A \cdot A) & = (\alpha_1, \ldots , \alpha_s) \cdot
\Delta \cdot (\alpha_1, \ldots , \alpha_s)^t  \\
 & = (\underline{w} \cdot P) \cdot \Delta \cdot (\underline{w}
 \cdot P)^t \\
 & =  -\underline{w} \cdot N \cdot \underline{w}^t  \\
 & = - \widehat{\underline{n}} \cdot M \cdot N \cdot
(\widehat{\underline{n}} \cdot
 M)^t  \\
 & = \widehat{\underline{n}} \cdot M^t \cdot
 \widehat{\underline{n}}^t \\
  & = \widehat{\underline{n}} \cdot M \cdot
  \widehat{\underline{n}}^t,
\end{align*}
as desired.
\qed
\medskip

\begin{Cor}
The codimension of the ideal $J^{D}(\underline{w})$ is
\[
h^{D}(\underline{w}) =
 \frac{1}{2} \left (
\sum_{i,i^{\prime}=1}^{s}m_{ii^{\prime}}\widehat{n}_i
\widehat{n}_{i^{\prime}} + \sum_{i=1}^{s} \widehat{n}_i \cdot
\left ( \sum_{{i^{\prime}}=1}^{s} m_{ii^{\prime}}
(2h_{i^{\prime}}-\overset{\bullet}{\nu_{i^{\prime}}})  - 1 \right
)
    \right ).
\]
\end{Cor}
\medskip

\begin{Not}
Let $I\subset I_0$, $J\subset J_0$, $\underline{n} \in
\mathcal{N}_{I,J}$ and $\underline{y}\in Y_{\underline{n}}$.
Recall that $\mathrm{\textsf{Init}}^{-1}(\underline{y})$ is an
affine space in $\mathbb{P}R^{\ast}$ of finite codimension
$F(\underline{n})$. Let $\underline{n} = (n_i, n'_{\sigma},
n''_{\sigma}, \tilde{n}'_j, \tilde{n}''_j )$ ; $1\le i\le s$,
$\sigma\in I$, $j\in J$. By considering $\widehat{\underline{n}} \in
\mathbb{Z}_{\ge 0}^s$ the element with entries
\[
\widehat{n}_i :=  n_i + \sum_{\substack{\sigma \in I
\\ i_1(\sigma)=i}} n_{\sigma}^{\prime} + \sum_{\substack{\sigma \in
I \\ i_2(\sigma)=i}} n_{\sigma}^{\prime \prime} +
\sum_{\substack{j \in J \\ i_1(j)=i }} \widetilde{n}_j^{\prime},
\]
we define
\begin{align*}
  \underline{w}(\underline{n}) & := \widehat{\underline{n}}
\cdot M  \\
  v_j (\underline{n}) & := w_{i_1(j)}(\underline{n}) +
\widetilde{n}_j^{\prime \prime} \cdot h_{i_1(j)}
\end{align*}
for $\underline{n} \in \mathcal{N}_{I,J}$ and $1 \le j \le r$.
Note that $v_j (\underline{n})$ is nothing but the order of a
function $g$ such that $\mathrm{\textsf{Init}}(g) \in
Y_{\underline{n}}$ on the component $C_j$ of the curve,
$w_{i}(\underline{n})$ is the multiplicity along the component
$E_{i}$ of the exceptional divisor of the lifting of a function
$g$ such that $\mathrm{\textsf{Init}}(g) \in Y_{\underline{n}}$,
and $\widetilde{n}^{\prime \prime}_j = \widetilde{n}^{\prime
\prime}_j (g)$.
\end{Not}
\medskip

\begin{Prop} \label{prop:efedeene}
We have
\begin{align*}
F(\underline{n}) & = \frac{1}{2} \left (
\sum_{i,i^{\prime}=1}^{s}m_{ii^{\prime}}\widehat{n}_i
\widehat{n}_{i^{\prime}} + \sum_{i=1}^{s} \widehat{n}_i \cdot
\left ( \sum_{{i^{\prime}}=1}^{s} m_{ii^{\prime}}
(2h_{i^{\prime}}-\overset{\bullet}{\nu_{i^{\prime}}}) +(2
h_{i^{\prime}} - 1) \right )
    \right ) \\
    &  + \sum_{j \in J} \widetilde{n}_j^{\prime \prime} \cdot
    h_j.
\end{align*}
\end{Prop}

\dem~The codimension
$F(\underline{n})$ is equal to the codimension in $R$ of the ideal
\[
I_{\underline{n}} = \left \{z \in R^{\ast} \mid w_i (z) > w_i
(\underline{n}), ~ 1 \le i \le s; ~v_j (z) > v_j (\underline{n}),~
1 \le j \le r \right \}
\]
minus $1$ by Corollary \ref{cor:F(n)}. Making additional blow--ups at the intersection
points of the strict transforms of the curve $C$ with the
exceptional divisor reduces our problem to the case $J
=\varnothing$. Then we have to compute the codimension $h^D
(\underline{w}(\underline{n})+\underline{1})-1$, with
\[
h^D (\underline{w}):= \dim_{k_R} R/J^D (\underline{w})
\]
for
$\underline{w} = (w_1, \ldots, w_s) \in \mathbb{Z}^s$ and $J^D(\underline{w}):= \left \{ z \in R^{\ast} \mid w_i (z) \ge
w_i, ~ ~ 1 \le i \le s \right \}$. By Proposition \ref{H-De},
Lemma \ref{lem:DK} and Lemma \ref{lem:DD}, we get
\begin{align*}
h^D (\underline{w}(\underline{n}))& = - \frac{1}{2} \left ( \deg_X
(A \cdot A) + \deg_X (A \cdot \mathbb{K}) \right )\\
& = \frac{1}{2} \left ( \widehat{\underline{n}} \cdot M \cdot
\widehat{\underline{n}}^t - \widehat{\underline{n}} \cdot
\underline{1}^t + \widehat{\underline{n}} \cdot M \cdot
\underline{\varepsilon}^t \right ).
\end{align*}
The above formula cannot be applied to compute $h^D(\underline{w}
+ \underline{1})$ because in general $\underline{w} +
\underline{1}$ is not of the form $\widehat{\underline{n}} \cdot M$ (i.e. does
not belong to the semigroup defined by the divisorial valuations
considered). Of course, in such a case the computation follows in
the same way. The map $\mathrm{\textsf{Init}}$ induces a fibration
\[
\mathbb{P} J^{D} (\underline{w} (\underline{n}))\supset Z
\longrightarrow \prod_{i=1}^{s} \mathcal{S}^{\widehat{n}_i}
\overset{\bullet}{E_i} = W,
\]
where $Z = \mathrm{\textsf{Init}}^{-1}(W)$ is an open subset of
$\mathbb{P}J^D(\underline{w}(\underline{n}))$. The fibre is the
set of elements in $Z$ whose image via the map
$\mathrm{\textsf{Init}}$ coincides (i.e., the map
$\mathrm{\textsf{Init}}$ applied to a non--zero element $g \in
\mathbb{P}R$ with $I(g)=\varnothing$); this set is by Proposition
\ref{lemma:clave} equal to $J^D (\underline{w}
(\underline{n})+\underline{1})$, and then
\[
h^D (\underline{w}(\underline{n})+\underline{1}) = 1 +  h^D
(\underline{w}(\underline{n})) + \dim_{k_R} \left (
\prod_{i=1}^{s} \mathcal{S}^{\widehat{n}_i} \overset{\bullet}{E_i}
\right );
\]
hence
\[
F(\underline{n})  = h^D (\underline{w}(\underline{n})) +
\dim_{k_R} \left ( \prod_{i=1}^{s} \mathcal{S}^{\widehat{n}_i}
\overset{\bullet}{E_i}
\right ) = h^D (\underline{w}(\underline{n})) + \sum_{i=1}^s
\widehat{n_i} h_i \;
\]
and the formula follows straightforward. \qed
\medskip

\subsection{Explicit formula.} \label{subsection:final}

\Nr \label{nr:2stars} The generalised Euler characteristic
satisfies the Fubini rule (see for instance \cite[\S3, pp.~
128--129]{viro}). Because of Corollary \ref{cor:F(n)} and Lemma
\ref{lem:36}, we can apply Fubini's formula to the map
\textsf{Init} in order to get
\[
P_g (t_1, \ldots ,t_r;\mathbb{L})= \int_{Y}
\mathbb{L}^{-F(\underline{n})}
\underline{t}^{\underline{v}(\underline{n})} d \chi_g = \sum_{\substack{I \subset I_0 \\
J \subset J_0}} ~ \sum_{\underline{n} \in \mathcal{N}_{I,J}}
\mathbb{L}^{-F(\underline{n})} [Y_{\underline{n}}] \cdot
\underline{t}^{\underline{v}(\underline{n})}. \eqno(**)
\]

By \ref{Yn} we have $Y_{\underline{n}} = \prod_{i=1}^{s} \mathcal{S}^{n_i}
\overset{\circ}{E_i} \times \prod_{\sigma \in I}
{k}_{\sigma}^{\ast} \times \prod_{j \in J} {k}_{j}^{\ast}$, 
and the class of $Y_{\underline{n}}$ in the Grothendieck ring is
\[
[Y_{\underline{n}}] = \prod_{i=1}^s [\mathcal{S}^{n_i}
\overset{\circ}{E_i}] \cdot \prod_{\sigma \in I}
[k_{\sigma}^{\ast}] \cdot \prod_{j \in J} [k_{j}^{\ast}],
\]
where $[k_{\sigma}^{\ast}]=[k_{\sigma}]-1 = [\mathrm{Spec}(k_{\sigma})]
\mathbb{L}-1$ and $[k_{j}^{\ast}]=[k_{j}]-1 = [\mathrm{Spec}(k_j)] \mathbb{L}-1$
(cf. \ref{nr:38}), and $[\mathcal{S}^{n_i} \overset{\circ}{E_i}]$
is computed in the following lemma.
\medskip

\begin{Lem}
Preserving notations as above, we have
\[
[\mathcal{S}^{n_i} \overset{\circ}{E_i}]= [\mathrm{Spec}(k_i)]^{n_i}
\mathbb{L}^{n_i} \sum_{l=0}^{\min \{n_i,
\overset{\circ}{\nu_i} -1 \}} (-1)^{l} [\mathrm{Spec}(k_i)]^{-l} {
\overset{\circ}{\nu_i} -1 \choose l } \mathbb{L}^{-l}.
\]
\end{Lem}

\dem~Proceeding as in 
\cite[Theorem 1, p.~51]{gzlume}, since $[\overset{\circ}{E_i}]= [\mathrm{Spec}(k_i)] \mathbb{L} + 1 - \overset{\circ}{\nu_i}$, we obtain
\begin{align*}
 \sum_{l=0}^{\infty} [\mathcal{S}^l \overset{\circ}{E_i}] t^l & = (1-t)^{-[\overset{\circ}{E_i}]} \\
 & = (1-t)^{- ( [\mathrm{Spec}(k_i)]\mathbb{L} + 1 - \overset{\circ}{\nu_i} )} \\
 & = (1-t)^{-[\mathrm{Spec}(k_i)]\mathbb{L}} (1-t)^{\overset{\circ}{\nu_i} - 1} \\
  & = \sum_{l=0}^{\infty} [\mathrm{Spec}(k_i)]^{l} \mathbb{L}^{l} t^l \cdot (1-t)^{\overset{\circ}{\nu_i}-1}.
\end{align*}
Therefore
\begin{align*}
[\mathcal{S}^{n_i} \overset{\circ}{E_i}] & =  \sum_{l=0}^{\min
\{n_i, \overset{\circ}{\nu_i} - 1 \}} (-1)^l
([\mathrm{Spec}(k_i)]\mathbb{L})^{n_i - l}
{\overset{\circ}{\nu_i} - 1 \choose l } \\
  & =  [\mathrm{Spec}(k_i)]^{n_i} \mathbb{L}^{n_i} \sum_{l=0}^{\min
\{n_i, \overset{\circ}{\nu_i} -1 \}} (-1)^{l} [\mathrm{Spec}(k_i)]^{-l}
{\overset{\circ}{\nu_i} -1 \choose l } \mathbb{L}^{-l} .
\end{align*}
\qed
\medskip

Once the classes in the Grothendieck ring have been computed, by
Equation ($\ast \ast$) in \ref{nr:2stars} we can conclude the
description of the generalised Poincar\'e series in terms of an
embedded resolution:
\medskip

\begin{Theo} \label{thm:uno}
\begin{eqnarray}
P_g(t_1, \ldots, t_r ; \mathbb{L}) & = & \sum_{\substack{I \subset I_0 \\
J \subset J_0}}~ \sum_{\underline{n} \in \mathcal{N}_{I,J}}
\mathbb{L}^{\sum_{i=1}^{s} n_i -F(\underline{n})}\cdot S_I
(\mathbb{L}) \cdot S_J (\mathbb{L}) \cdot \prod_{i=1}^{s}[\mathrm{Spec}(k_i)]^{n_i}
\nonumber \\
  & \times & \prod_{i=1}^s \left ( \sum_{l=0}^{\min
\{n_i, \overset{\circ}{\nu_i} -1 \}} (-1)^{l} [\mathrm{Spec}(k_i)]^{-l}
{\overset{\circ}{\nu_i} -1 \choose l } \mathbb{L}^{-l} \right
) \underline{t}^{\underline{v}(\underline{n})}, \nonumber
\end{eqnarray}
where $S_I(\mathbb{L}):=\prod_{\sigma \in I} ([\mathrm{Spec}(k_{\sigma})]
\mathbb{L}-1)$ and $S_J (\mathbb{L}):= \prod_{j \in J} ([\mathrm{Spec}(k_j)]
\mathbb{L}-1)$.
\end{Theo}
\medskip

Notice that $F(\underline{n})$ has been already computed in
Proposition \ref{prop:efedeene} in terms of the intersection
matrix. This result generalises Theorem 1 of \cite{cadegu11}, which
holds only for complex curve singularities. Indeed, the formula
turns out to be easier for the \emph{totally rational} case, i.e. if all field extensions have degree one (in other words,
$h_{\sigma}=h_j=1$ for all $\sigma \in I$ and $j \in J$):
\medskip

\begin{Cor}[Campillo, Delgado, Gusein--Zade]
Assume the ring $R$ to be totally rational. We have
\begin{eqnarray}
P_g(t_1, \ldots, t_r ; \mathbb{L}) & = & \sum_{\substack{I \subset I_0 \\
J \subset J_0}} ~ \sum_{\underline{n} \in \mathcal{N}_{I,J}}
\mathbb{L}^{\sharp (I) + \sharp (J) + \sum_{i=1}^{s} n_i
-F(\underline{n})} \left (1-\mathbb{L}^{-1} \right )^{\sharp (I) +
\sharp (J)}
\nonumber \\
  & \times & \prod_{i=1}^s \left ( \sum_{l=0}^{\min
\{n_i, \overset{\circ}{\nu_i} -1 \}} (-1)^{l}
{\overset{\circ}{\nu_i} -1 \choose l } \mathbb{L}^{-l} \right )
\underline{t}^{\underline{v}(\underline{n})}. \nonumber
\end{eqnarray}
\end{Cor}
\medskip

\begin{Rem}
The specialisation $\mathbb{L} \rightarrow 1$ gives the connection
with the classical Poincar\'e series $P(\underline{t})$:
\[
P_g (t_1, \ldots ,t_r; 1) = \sum_{\substack{I \subset I_0 \\
J \subset J_0}} ~ ~ ~ \sum_{\underline{n} \in \mathcal{N}_{I,J}}
\underline{t}^{\underline{v}(\underline{n})} =  \int_Y
\underline{t}^{\underline{v}(\underline{n})} d \chi = P
(\underline{t}).
\]
\end{Rem}
\medskip

\section{Divisorial Poincar\'e series and resolution} \label{sec:div}

Let assume the ring $R$ to have a perfect coefficient field $K$
along this section.
\medskip

\begin{Defi}
A divisorial valuation of $R$ is a discrete valuation of rank one
(with group of values $\mathbb{Z}$) of the field of fractions of
$R$ with $R \cap \mm_{\nu}=\mm$, if $R_{\nu}$ is the valuation
ring with maximal ideal $\mm_{\nu}$, and with transcendence degree
of the residual extension $R_{\nu}/\mm_{\nu} : R / \mm$ equal to
$1$.
\end{Defi}
\medskip

\begin{Not} \label{not:52}
Let $D:=\{w_1 , \ldots , w_s \}$ be a finite set of divisorial
valuations of $R$, and let $W_i$ denote the discrete valuation ring associated with $w_i$ for all $1 \le i \le s$. 
The divisorial value semigroup associated with
$D$ is the additive sub--semigroup $S^D$ of $\mathbb{Z}_{\ge 0}^s$
given by
\[
S^D := \{ \underline{w}(z):= (w_1 (z), \ldots , w_s (z)) \mid z\in
R \setminus \{ 0 \} \}.
\]

For $I \subset I_0$ consider the set
\begin{align*}
\mathcal{N}_{I}^D& :=  \{ \underline{n}:= (n_i,
n_{\sigma}^{\prime}, n_{\sigma}^{\prime \prime}) \mid n_i \ge 0,~1
\le i \le s;  ~ n^{\prime}_{\sigma}>0, ~ n^{\prime
\prime}_{\sigma}>0,~ \sigma \in I \};
\end{align*}
(the super--index $D$ of $\mathcal{N}_I^{D}$ and other next
notations refers to the word ``divisorial''). Consider the ideals
$J^D (\underline{w}):= \{z \in R \mid \underline{w}(z) \ge
\underline{w} \}$ already defined in \ref{nr:413}. We introduce a
definition similar to \ref{poincareseries} in the divisorial
context (cf. \cite{Delgado5}):
\end{Not}
\medskip

\begin{Defi}  \label{defn:generaliseddivisorial}
The generalized divisorial Poincar\'e series of the multi--index
filtration given by the ideals $J^D(\underline{w})$ is defined to
be the integral
\begin{equation} \label{eqn:poincareseriesdef}
P^D_g (t_1, \ldots ,t_s;\mathbb{L}):=\int_{\mathbb{P} R}
\underline{t}^{\underline{w}(z)} d \chi_g \in \mathcal{M}_{k_R}.
\nonumber
\end{equation}
\end{Defi}
\medskip

We can express this divisorial series in terms of an embedded
resolution of curves with a similar argument as that used in
Theorem \ref{thm:uno}:
\medskip

\begin{Theo} \label{thm:gendiv1}
\begin{eqnarray}
P^D_g(t_1, \ldots, t_s ; \mathbb{L}) & = & \sum_{\substack{I \subset I_0 \\
J \subset J_0}} ~ ~ ~ \sum_{\underline{n} \in \mathcal{N}_{I}}
\mathbb{L}^{\sum_{i=1}^{s} n_i -F^D(\underline{n})}\cdot S_I
(\mathbb{L}) \cdot \prod_{i=1}^{s}[\mathrm{Spec}(k_i)]^{n_i} \times
\nonumber \\
  & \times & \prod_{i=1}^s \left ( \sum_{l=0}^{\min
\{n_i, \overset{\bullet}{\nu_i} -1 \}} (-1)^{l} [\mathrm{Spec}(k_i)]^{-l}
{\overset{\bullet}{\nu_i} -1 \choose l } \mathbb{L}^{-l}
\right ) \underline{t}^{\underline{w}(\underline{n})}, \nonumber
\end{eqnarray}
where $S_I(\mathbb{L}):=\prod_{\sigma \in I} ([\mathrm{Spec}(k_{\sigma})]
\mathbb{L}-1)$ and, for $\underline{n} \in \mathcal{N}^D_I$ we denote
\begin{itemize}
    \item $ \widehat{n}_i :=  n_i + \sum_{\substack{\sigma \in I
\\ i_1(\sigma)=i}} n_{\sigma}^{\prime} + \sum_{\substack{\sigma \in
I \\ i_2(\sigma)=i}} n_{\sigma}^{\prime \prime}$
    \item $\underline{w}(\underline{n}):= \widehat{\underline{n}} \cdot M \nonumber$
    \end{itemize}
and the codimension $F^D(\underline{n})$ is equal to
\[
\frac{1}{2} \left (
\sum_{i,i^{\prime}=1}^{s}m_{ii^{\prime}}\widehat{n}_i
\widehat{n}_{i^{\prime}} + \sum_{i=1}^{s} \widehat{n}_i \cdot
\left ( \sum_{{i^{\prime}}=1}^{s} m_{ii^{\prime}}
(2h_{i^{\prime}}-\overset{\bullet}{\nu_{i^{\prime}}}) +(2
h_{i^{\prime}} - 1) \right )
    \right ).
    \]
\end{Theo}

\dem~We proceed as in Theorem \ref{thm:uno}, but taking the space
\[
Y^{D}:=\bigcup_{I \subset I_0} \bigcup_{\underline{n} \in
\mathcal{N}_I^{D}}  Y_{\underline{n}},
\]
with $Y_{\underline{n}}:= \prod_{i=1}^{s} \mathcal{S}^{n_i}
\overset{\bullet}{E_i} \times \prod_{\sigma \in I}
k^{\ast}_{\sigma}$, instead of the space $Y$ defined in Equation (\dag) of $\ref{Yn}$. \qed

\Nr Let us now consider the extended semigroup $\widehat{S}_D$
coming from the divisorial filtration $\{J^D(\underline{w})
\}$ (cf. \cite{Delgado5}) and the Poincar\'e series
\[
\widehat{P}_g^{D}(t_1, \ldots ,t_s;\mathbb{L}):=
\int_{\widehat{S}_D} \underline{t}^{\underline{v}(h)} d \chi_g.
\]
\Nr \label{nr:56}  The next goal is to express this Poincar\'e series in terms of an
embedded resolution. First, we define a semigroup homomorphism $
\Pi: Y^{D} \rightarrow \pro \widehat{S}_D$ 
in the following way: For $y \in Y^D$, $y$ is represented by a set
of smooth closed points of $\overset{\bullet}{E}=\bigcup_{i=1}^{s}
\overset{\bullet}{E_i}$ with $n_i$ points $Q_1^{i}, \ldots ,
Q_{n_i}^{i}$ on the component $\overset{\bullet}{E_i}$. Let $Q^{i}$ be a point of $\overset{\bullet}{E_i}$. Assume
$\overset{\bullet}{E_i}$ to have a local equation $\{x_i=0 \}$ at
$Q^{i}$. Let us take an element $\widetilde{\gamma}_{Q^{i}}$
meeting $E_i$ at $Q^{i}$ transversally, and let
$R_{Q^{i}}:=\oo_{X,Q^{i}}$. By \cite[Lemma 3.1.11]{moyano} (see
also \cite[Section 6]{moyano2}), there exists an element
$g_{Q^{i}} \in R$ such that its strict transform at $Q^{i}$ is
$\widetilde{\gamma}_{Q^{i}}$. Recall now the existence of a the one--to--one correspondence
between the components of the exceptional divisor and the set $D$ (see \cite[Chap. VII]{kiyek}).
If we consider the transform of
$g_{Q^{i}}$ in the discrete valuation ring $W_i$ corresponding to
the component $E_i$ (cf. \ref{not:52}), let us say
$\widetilde{g}_{Q^{i}} \in W_i$, then there exists an uniformising
parameter $t_i \in W_i$ such that $\widetilde{g}_{Q^{i}}=\alpha_i
(g_{Q^{i}}) \cdot t_i^{z_i}$, where $\alpha_i(g_{Q^{i}})$ is a
unit of $W_i$ and $z_i=w_i (g_{Q^{i}})$, for all $1 \le i \le s$.
If we take the image $a_i (g_{Q^{i}})$ of $\alpha_i (g_{Q^{i}})$
in $k_i^{\ast}$, then the element $(a_i (g_{Q^{i}}), w_i
(g_{Q^{i}}))$ belongs to the divisorial semigroup $\widehat{S}_D$.
\medskip

By definition, $\Pi(y) \in \pro \widehat{S}_D$ is represented by
the element $(\underline{v}(g),\underline{a}(g)) \in
\widehat{S}_D$, where $g= \prod_{i=1}^{s} \prod_{j=1}^{n_i}
g_{Q_j^{i}}$. But it is independent of the chosen representant
$\widetilde{\gamma}_{Q^{i}}$ at a point $Q^{i}$ of the component
$E_i$ of the exceptional divisor, as the following result shows.
\medskip

\begin{Lem} \label{lem:57}
The element $\Pi(y)$ belonging to the projectivisation  $\pro
\widehat{S}_D$ of the divisorial extended semigroup
$\widehat{S}_D$ does not depend on the choice of the curves
$\widetilde{\gamma}_{P}$, with $P$ a point at
$\overset{\bullet}{E_i}$ for any $i \in \{1, \ldots , s \}$.
\end{Lem}

\dem~Let $\widetilde{\gamma}^{\prime}_{P}$ be a transversal
element to $E$ (a curvette) at a closed point $P \in
\overset{\bullet}{E_i}$ for some $i \in \{1, \ldots , s \}$ coming
from an element $g^{\prime}_{P} \in R$ and let
$g^{\prime}=\prod_{i= 1}^{s} \prod_{j=1}^{n_i}
g^{\prime}_{Q_j^{i}}$. Let $\widetilde{g}=g \circ \pi$ and
$\widetilde{g^{\prime}}=g^{\prime} \circ \pi$ be the liftings of
the functions $g$ and $g^{\prime}$ to the space $X$ of the
resolution. We can set the function
$\psi=\frac{\widetilde{g^{\prime}}}{\widetilde{g}}$, whose
restriction to every component $E_i$ of $E$ is regular, and
therefore constant (since the $E_i$ are projective lines).
A similar argument as that in the proof of Proposition
\ref{lemma:clave} shows that $\psi$ is constant on $E$ and equal
to some $\alpha \in k_R$. It implies that
$\underline{v}(g^{\prime})=\underline{v}(\alpha
g)=\underline{v}(g)$ and $\underline{a}(g^{\prime})=\alpha \cdot
\underline{a}(g)$. \qed
\medskip

\begin{Lem} \label{lemma:iso}
The semigroup homomorphism $\Pi$ is an isomorphism.
\end{Lem}
\dem~It is a straight consequence of Proposition
\ref{lemma:clave}. \qed
\medskip

\begin{Theo} \label{thm:divisorialintermsof}
If we denote $\underline{t}^{\underline{m}_{\sigma}}:=
(1-\underline{t}^{\underline{m}_{i_1(\sigma)}})(1-\underline{t}^{\underline{m}_{i_2(\sigma)}})$,
then the following equalities hold:
\begin{align*}
\widehat{P}_g^{D}(t_1, \ldots ,t_s;\mathbb{L}) & = \chi_g (\pro
\widehat{S}_D)\\
& = \int_{\pro \widehat{S}_D}
\underline{t}^{\underline{w}}d \chi_g \\
& =  \frac{\prod_{\sigma \in I_0} \left (
\underline{t}^{\underline{m}_{\sigma}} \right )^{h_{\sigma}} +
\left ( \underline{t}^{\underline{m}_{\sigma}} \right
)^{h_{\sigma}-1} ([\mathrm{Spec}(k_{\sigma})] \mathbb{L}-1)
\underline{t}^{\underline{m}_{i_1(\sigma)}}
\underline{t}^{\underline{m}_{i_2(\sigma)}} }{\prod_{i=1}^{s}
(1- \underline{t}^{\underline{m}_i})(1- [\mathrm{Spec}(k_{i})]  \mathbb{L}
\underline{t}^{\underline{m}_i})}.
\end{align*}
\end{Theo}

\dem~By Lemma \ref{lemma:iso}, we have
\begin{align*}
\int_{\pro \widehat{S}_D} \underline{t}^{\underline{w}} d \chi_g   & = \int_{Y^D} \underline{t}^{\underline{w}} d \chi_g =   \sum_{I \subset I_0} \sum_{\underline{n} \in \mathcal{N}_I^D} [Y^D_{\underline{n}}] \cdot \underline{t}^{\underline{w}(\underline{n})}  \\
   & = \sum_{I \subset I_0} \sum_{\underline{n} \in \mathcal{N}_I^D} \prod_{\sigma \in I} ([\mathrm{Spec}(k_{\sigma})] \mathbb{L}-1) \prod_{i=1}^{s}[\mathcal{S}^{n_i}\overset{\bullet}{E_i}] \cdot \underline{t}^{\sum_{i=1}^{s} \widehat{n}_i \underline{m}_i}\\
   & = A(\underline{t}) \cdot B(\underline{t}),
\end{align*}   
where
\begin{align*}
A(\underline{t}) & :=\sum_{\substack{n_i \ge 0 \\ i=1,\ldots,s}} \prod_{i=1}^{s}[\mathcal{S}^{n_i} \overset{\bullet}{E_i}] \cdot \underline{t}^{\sum_{i=1}^{s}n_i \underline{m}_i} \\
B(\underline{t}) & := \sum_{I \subset I_0} \prod_{\sigma \in I} ([\mathrm{Spec}(k_{\sigma})]\mathbb{L}-1) \sum_{\substack{n^{\prime}_{\sigma}>0~ n^{\prime \prime}_{\sigma}>0 \\ \sigma \in I}} \underline{t}^{n^{\prime}_{\sigma} \underline{m}_{i_1(\sigma)}} \cdot \underline{t}^{n^{\prime \prime}_{\sigma} \cdot \underline{m}_{i_2(\sigma)}}. 
 \end{align*}
 
Concerning $A(\underline{t})$ we refer to
\cite[Theorem 1]{gzlume} to get
\begin{align*}
A (\underline{t})   & = \prod_{i=1}^{s} \left ( \sum_{n=0}^{\infty} [\mathcal{S}^{n} \overset{\bullet}{E_i}] \cdot \underline{t}^{n \cdot \underline{m}_i} \right )=\prod_{i=1}^{s} (1-\underline{t}^{\underline{m}_i})^{-[\overset{\bullet}{E_i}]}\\
   & = \prod_{i=1}^{s} (1-\underline{t}^{\underline{m}_i})^{-[\pro^{1}_{k_i}]}
   \cdot(1-\underline{t}^{\underline{m}_i})^{\overset{\bullet}{\nu_i}}.
\end{align*}

Since $[\pro^{1}_{k_i}] = [\mathrm{Spec}(k_{i})] \mathbb{L}+1$ and
$(1-\underline{t})^{- \mathbb{L}} = (1-
\mathbb{L}\underline{t})^{-1}$ (cf. \cite{gzlume}), we obtain
\begin{align*}
 \left ( 1-\underline{t}^{\underline{m}_i} \right
)^{-[\pro^{1}_{k_i}]} & = \left ( 1 -
\underline{t}^{\underline{m}_i} \right )^{-[\mathrm{Spec}(k_{i})]  \mathbb{L}} \cdot
\left ( 1-\underline{t}^{\underline{m}_i} \right
)^{-1} \\
& = (1- [\mathrm{Spec}(k_{i})] \mathbb{L}\underline{t}^{\underline{m}_i})^{-1}
\cdot (1-\underline{t}^{\underline{m}_i})^{-1}.
\end{align*}
On the other hand, by the equality
\[
(1-\underline{t}^{\underline{m}_i})^{\overset{\bullet}{\nu_i}} =
\prod_{\sigma \in I_0}  (1-
\underline{t}^{\underline{m}_{i_1(\sigma)}} )^{h_{\sigma}} (1-
\underline{t}^{\underline{m}_{i_2(\sigma)}})^{h_{\sigma}} ,
\]
one deduces
\begin{align*}
 A(\underline{t}) & = \prod_{i=1}^{s} \frac{1}{\left ( 1- \underline{t}^{\underline{m}_i} \right ) \left ( 1- [\mathrm{Spec}(k_{i})] \mathbb{L} \underline{t}^{\underline{m}_i} \right
 )} \cdot \prod_{\sigma \in I_0} \left ( \left ( 1- \underline{t}^{\underline{m}_{i_1(\sigma)}} \right ) \left ( 1- \underline{t}^{\underline{m}_{i_2(\sigma)}} \right
 ) \right )^{h_{\sigma}}.
\end{align*}

The factor $B(\underline{t}) $ is equal to
\begin{align*}
 B(\underline{t})  & =  \sum_{I \subset I_0} \prod_{\sigma \in I} ([\mathrm{Spec}(k_{\sigma})] \mathbb{L}-1) \prod_{\sigma \in I} \frac{\underline{t}^{\underline{m}_{i_1(\sigma)}}}{1-\underline{t}^{\underline{m}_{i_1(\sigma)}}} \cdot \frac{\underline{t}^{\underline{m}_{i_2(\sigma)}}}{1-\underline{t}^{\underline{m}_{i_2(\sigma)}}} \\
 & =  \prod_{\sigma \in I_0} \left ( 1 + ([\mathrm{Spec}(k_{\sigma})] \mathbb{L}-1) \frac{\underline{t}^{\underline{m}_{i_1(\sigma)}}}{1-\underline{t}^{\underline{m}_{i_1(\sigma)}}} \cdot \frac{\underline{t}^{\underline{m}_{i_2(\sigma)}}}{1-\underline{t}^{\underline{m}_{i_2(\sigma)}}} \right )\\
 & =  \prod_{\sigma \in I_0} \frac{(1-\underline{t}^{\underline{m}_{i_1(\sigma)}})(1-\underline{t}^{\underline{m}_{i_2(\sigma)}}) + ([\mathrm{Spec}(k_{\sigma})]  \mathbb{L}-1)\underline{t}^{\underline{m}_{i_1(\sigma)}} \underline{t}^{\underline{m}_{i_2(\sigma)}}
 }{(1-\underline{t}^{\underline{m}_{i_1(\sigma)}})(1-\underline{t}^{\underline{m}_{i_2(\sigma)}})}.
\end{align*}
Denoting
$\underline{t}^{\underline{m}_{\sigma}}:=(1-\underline{t}^{\underline{m}_{i_1(\sigma)}})(1-\underline{t}^{\underline{m}_{i_2(\sigma)}})$,
we get
\begin{align*}
A(\underline{t})  \cdot B(\underline{t})  & =  \frac{\prod_{\sigma \in I_0} \left (
\underline{t}^{\underline{m}_{\sigma}} \right )^{h_{\sigma}-1}+
\left ( \underline{t}^{\underline{m}_{\sigma}} + ([\mathrm{Spec}(k_{\sigma})] 
\mathbb{L}-1) \underline{t}^{\underline{m}_{i_1(\sigma)}}
\underline{t}^{\underline{m}_{i_2(\sigma)}} \right
)}{\prod_{i=1}^{s}  (1- \underline{t}^{\underline{m}_i})(1-[\mathrm{Spec}(k_{i})] 
\mathbb{L} \underline{t}^{\underline{m}_i})}\\
& = \frac{\prod_{\sigma \in I_0} \left (
\underline{t}^{\underline{m}_{\sigma}} \right )^{h_{\sigma}} +
\left ( \underline{t}^{\underline{m}_{\sigma}} \right
)^{h_{\sigma}-1} ([\mathrm{Spec}(k_{\sigma})]  \mathbb{L}-1)
\underline{t}^{\underline{m}_{i_1(\sigma)}}
\underline{t}^{\underline{m}_{i_2(\sigma)}} }{\prod_{i=1}^{s}
(1- \underline{t}^{\underline{m}_i})(1- [\mathrm{Spec}(k_{i})] \mathbb{L}
\underline{t}^{\underline{m}_i})}.
\end{align*}
\qed
\medskip

\begin{Cor}[Campillo, Delgado, Gusein--Zade]
If the ring $R$ is totally rational, then we have
$$
\widehat{P}_g^{D}(t_1, \ldots ,t_s;\mathbb{L}) =
\frac{\prod_{\sigma \in I_0} \left (
1-\underline{t}^{\underline{m}_{i_1(\sigma)}}
-\underline{t}^{\underline{m}_{i_2(\sigma)}} + \mathbb{L}
\underline{t}^{\underline{m}_{i_1(\sigma)}}
\underline{t}^{\underline{m}_{i_2(\sigma)}} \right
)}{\prod_{i=1}^{s}
(1-\underline{t}^{\underline{m}_i})(1-\mathbb{L}\underline{t}^{\underline{m}_i})}.
$$
\end{Cor}


\medskip



\begin{thebibliography}{99}

\bibitem{blickle} M.~Blickle: \emph{A short course in geometric motivic
integration}. Preprint arXiv: math.AG/0507404,~ 42 pp.~ (2005).

\bibitem{cadegu1} A.~Campillo, F.~ Delgado, S.~M.~Gusein--Zade: {\em
On the monodromy of a plane curve singularity and the Poincar\'e
series of its rings of functions}. Functional Analisis and its
Applications \textbf{33} (1),~56--57~(1999).

\bibitem{cadegu3} A.~Campillo, F.~ Delgado, S.~M.~Gusein--Zade: {\em
The Alexander polynomial of a plane curve singularity and the ring
of functions on it}. Russian Math. Surveys \textbf{54}(3),~634--635~(1999).

\bibitem{cadegu4} A.~Campillo, F.~ Delgado, S.~M.~Gusein--Zade: 
{\em Integration with res\-pect to the Euler characteristic over a
function space, and the Alexander polynomial of a plane curve
singularity}. Russian
Math. Surveys \textbf{55}(6),~ 1148--1149~(2000).

\bibitem{cadegu6} A.~Campillo, F.~ Delgado, S.~M.~Gusein--Zade: 
\emph{The Alexander polynomial of a plane curve singularity and
integrals with respect to the Euler characteristic}. Int.
Journal of Math. \textbf{14}(1),~ 47--54~ (2003).

\bibitem{duke} A.~Campillo, F.~ Delgado, S.~M.~Gusein--Zade: 
\emph{The Alexander polynomial of a plane curve singularity via
the ring of functions on it}. Duke Math. Journal \textbf{117} (1),~ 125--156~(2003).

\bibitem{cadegu11} A.~Campillo, F.~ Delgado, S.~M.~Gusein--Zade: 
\emph{Multi--index filtrations and motivic Poincar\'e series}.
Monatshefte f\"ur Mathematik \textbf{150}, ~193--209~ (2007).

\bibitem{cadeki} A.~Campillo, F. ~ Delgado, K.~Kiyek:  \emph{Gorenstein
 property and symmetry for one-dimensional local Cohen-Macaulay
 rings}. Manuscripta Math., \textbf{83}, ~ 405--423~(1994).

\bibitem{Delgado5} F.~Delgado de la Mata, S.~M.~Gusein--Zade: {\em
Poincar\'e series for several plane divisorial valuations}. Proc.
Edinburgh Math. Soc. \textbf{46}(2),~501--509~(2003).

\bibitem{Delgado4} F.~Delgado, C.~Galindo, A.~N\'u\~nez:
\emph{Generating sequences and Poincar\'e series for a finite set
of plane divisorial valuations}. Adv. Math. \textbf{219},~1632--1655~(2008).

\bibitem{demo} F.~Delgado de la Mata, J.~J.~Moyano--Fern\'andez: \emph{On the relation between the generalized Poincar\'e series and the St\"ohr Zeta
function}. Proc. Amer. Math. Soc. \textbf{137}(1),~51--59~(2009).

\bibitem{denef} J.~Denef, F.~Loeser:
\emph{Germs of arcs on singular algebraic varieties and motivic
integration}. Invent. Math. \textbf{135},~201--232~(1999).

\bibitem{hartshorne} R.~Hartshorne: \emph{Algebraic Geometry}.
Springer, New York ~(1977).

\bibitem{kimo} K.~Kiyek, J.~J.~Moyano--Fern\'andez: \emph{The Poincar\'{e}
series of a simple complete ideal of a two-dimensional regular
local ring}. J. Pure Appl. Algebra \textbf{213}(9),~1777--1787~(2009)

\bibitem{kiyek} K.~Kiyek, J.~L.~Vicente: \emph{Resolution of Curve and Surface Singularities in Characteristic
Zero}. Kluwer, Dordrecht~(2004).

\bibitem{ega1} A.~Grothendieck, J.~A.~Dieudonn\'e: \emph{\'El\'ements de G\'eom\'etrie Alg\'ebrique
 I.} Die Grundlehren der mathematischen Wissenschaften in
 Einzeldarstellungen. Band \textbf{166}. Springer, Berlin-Heidelberg-New
 York ~(1971).

\bibitem{gzlume} S.~M.~Gusein--Zade, I.~Luengo, A.~Melle:
\emph{A power structure over the Grothendieck ring of varieties}.
Mathematical Research Letters \textbf{11},~49--57~(2004).

\bibitem{liu} Q.~Liu: \emph{Algebraic Geometry and Arithmetic
Curves}. Oxford U.P., Oxford~(2002).

\bibitem{moyano} J.~J.~Moyano--Fern\'andez: \emph{Poincar\'e series associated with curves defined over finite
fields}. Dissertation, Universidad de Valladolid (Spain)~(2008).

\bibitem{moyano2} J.~J.~Moyano--Fern\'andez: \emph{Curvettes and
clusters of infinitely near points}. To appear in Rev. Mat. Complut. DOI:~10.1007/s13163-010-0048-1.

\bibitem{mozu} J.~J.~Moyano--Fern\'andez, W.~A.~Z\'u\~niga--Galindo: \emph{Motivic zeta functions for curve singularities}.  Nagoya Math. Journal \textbf{198},~47--75 ~(2010).

\bibitem{poonen} B.~Poonen: \emph{The Grothendieck ring of varieties is not a
domain}. Math. Res. Lett. \textbf{9}(4),~493--497~(2002).

\bibitem{viro} O.~Y.~Viro: \emph{Some integral calculus based on Euler
characteristic}. In: Viro, O.Y. (ed.) \emph{Topology and
Geometry--Rohlin Seminar.} Lect. Notes in Math. \textbf{1346},~pp.~127--138. Springer, Berlin, Heidelberg~(1988).
\end{thebibliography}
\end{document}